\documentclass[11pt]{article}
\usepackage[a4paper,left=1.5cm,right=1.5cm,top=2.5cm,bottom=2.5cm]{geometry}

\usepackage{amsmath}
\usepackage{amssymb}
\usepackage{amsfonts}
\usepackage{amsthm}
\usepackage{graphicx}
\usepackage{subcaption}
\usepackage{booktabs}
\usepackage{caption}
\usepackage{enumitem}
\usepackage{xr-hyper}
\usepackage[hidelinks]{hyperref}
\usepackage[nameinlink,capitalise,noabbrev]{cleveref}
\usepackage{amsopn}
\DeclareUnicodeCharacter{0301}{}

\DeclareGraphicsExtensions{.pdf,.png,.jpg,.eps}

\setlist[enumerate]{leftmargin=.5in}
\setlist[itemize]{leftmargin=.5in}

\newcommand{\headers}[2]{}
\newcommand{\email}[1]{\texttt{#1}}
\newenvironment{keywords}{\par\noindent\textbf{Keywords: }\ignorespaces}{\par}

\crefname{hypothesis}{Hypothesis}{Hypotheses}

\crefname{fact}{Fact}{Facts}

\title{Regular and chaotic Welander oscillations in a four-dimensional conceptual model for the Atlantic Meridional Overturning Circulation}

\author{John Bailie\thanks{Email: \email{john.bailie@auckland.ac.nz}.}
  \and Priya Subramanian
  \and Bernd Krauskopf}
\date{Department of Mathematics, University of Auckland, 
  Private Bag 92019, \\ Auckland 1142, New Zealand \\[5mm]\today
}

\begin{document}

\maketitle

\begin{abstract}
  The Atlantic Meridional Overturning Circulation (AMOC) is a key component of the Earth’s climate. Evidence indicates a twentieth-century weakening, and enhanced freshwater input to the subpolar North Atlantic may further reduce overturning strength. We present and study a conceptual four-dimensional, single-hemisphere box model with three compartments: a tropical surface box, a subpolar surface box, and a large deep-water box. Advective exchange couples the surface boxes and vertical exchange with the deep ocean is represented by a smooth convective-adjustment scheme. A comprehensive bifurcation analysis reveals an equilibrium structure with up to four coexisting overturning states, together with regimes of bistability and tristability. We identify families of periodic solutions and chaotic attractors with a clear timescale separation: a millennial oscillation is modulated by faster decadal-to-centennial variability arising from episodic shutdowns of subpolar convection. As prescribed freshwater fluxes increase, shutdown events become more frequent and the background overturning weakens. Additionally, for certain values of freshwater influx, the dynamics become chaotic, producing an irregular on--off switching of convection.
\end{abstract}






\section{Introduction}
The Atlantic Meridional Overturning Circulation (AMOC) is a large-scale ocean circulation spanning the entire Atlantic Ocean \cite{kuhlbrodt2007driving, lozier2010deconstructing, buckley2016observations}. Its upper limb transports warm, salty water from the tropics to the subpolar North Atlantic, particularly the Labrador, Irminger, and Nordic seas \cite{pickart2007impact}. There, surface cooling, salt rejection during sea-ice formation, and continued northward salt advection precondition the water column for deep convection to form the North Atlantic Deep Water (NADW) \cite{yashayaev2024intensification, marshall1999open}. A deep cold current brings the NADW  toward lower latitudes, where it returns to the upper ocean through wind-driven Ekman pumping and widespread diapycnal mixing (turbulent mixing that transports heat and salt across constant-density surfaces)~\cite{marshall2012closure, talley2013closure}. The upwelled waters rejoin the northward surface flow, completing the overturning circulation.

The AMOC is widely recognized as a critical tipping element of the climate system, whose substantial weakening or shutdown could trigger abrupt changes in other components, including the possibility of a tipping cascade in which additional climate subsystems are destabilized in turn~\cite{lenton2008tipping, wunderling2024climate}. Such a cascade could substantially disrupt the El Ni\~no--Southern Oscillation patterns~\cite{molina2022response, williamson2018effect} and lead to a weakened Gulf Stream that would reduce northward ocean heat transport and thereby cool the North Atlantic and northwestern Europe~\cite{vellinga2008impacts}. Moreover, reduced overturning is expected to weaken Atlantic nutrient transport and upwelling, leading to major impacts on marine ecosystems~\cite{lynch2024diminished}.

\begin{figure}[t!]
  \centering    \includegraphics[width=0.8\linewidth]{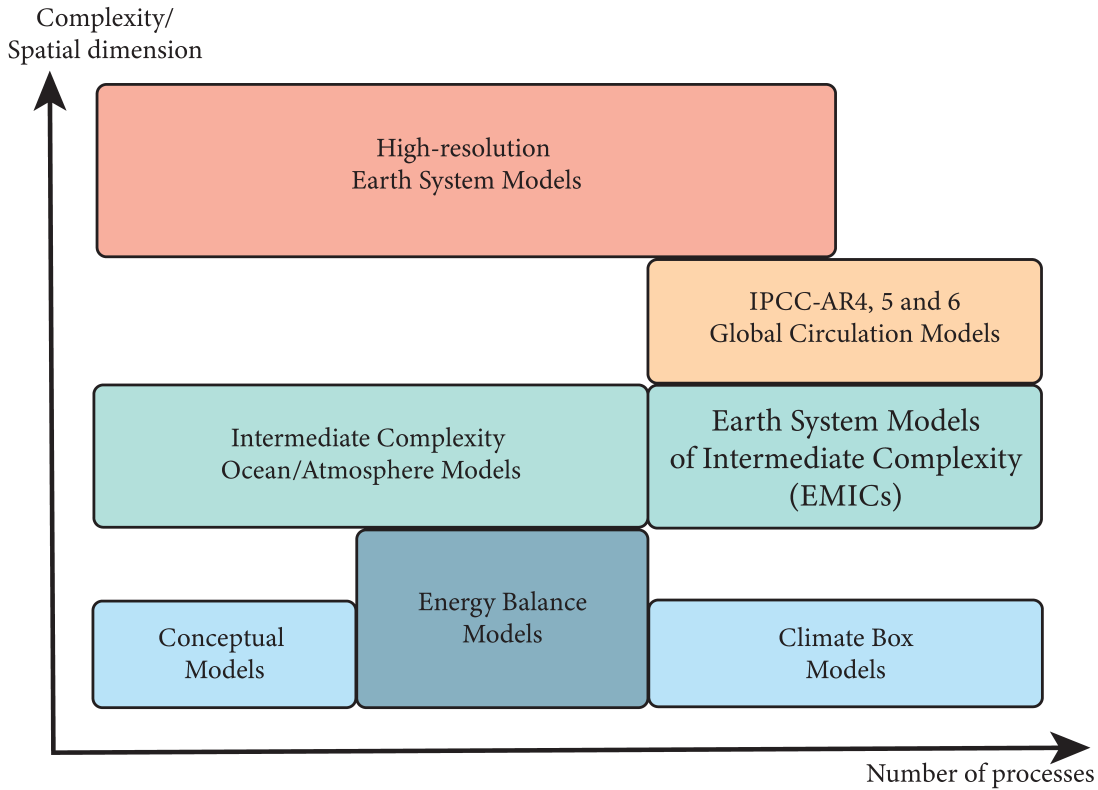}
  \caption{Hierarchy of climate models, organized by spatial dimensionality on the y-axis and by the number of climate processes on the x-axis. Adapted from~\cite{dijkstra2013nonlinear, dijkstra2024role}. }
  \label{fig:climate_heirarcy}
\end{figure}

Given these potentially severe impacts on the global climate, understanding how and when the AMOC might weaken requires careful modeling of the climate system. Figure~\ref{fig:climate_heirarcy} summarizes the hierarchy of models that climate scientists use to study the AMOC~\cite{dijkstra2013nonlinear, van2024role}. At the top, General Circulation Models (GCMs), typically run at horizontal resolutions of order $1^\circ$, simulate a fully coupled three-dimensional ocean--atmosphere system~\cite{boucher2020presentation, gent2011community, held2019structure, kuhlbrodt2018low}, while Earth System Models (ESMs) additionally incorporate biogeochemical processes~\cite{danabasoglu2020community, mauritsen2019developments, sellar2019ukesm1}. Their comprehensive scope enables a relatively faithful description of the global climate, yet can turn them into \emph{black boxes} that obscure individual mechanisms and limit the direct application of dynamical systems analysis~\cite{kuznetsov1998elements, meiss2007differential, wiggins2003introduction}. Earth system models of intermediate complexity (EMICs) provide coarser representations, typically at horizontal resolutions between $3^\circ$ and $7.5^\circ$, at a lower computational cost~\cite{goosse2010description, petoukhov2000climber, weaver2001uvic}; this is achieved, for instance, by reducing spatial and temporal resolution or by averaging zonally dependent processes~\cite{stocker1992zonally, hovine1994zonally}. They still retain many interacting components, which can make isolating specific mechanisms challenging. At the lower end of the modeling hierarchy lie conceptual box models: phenomenological descriptions that isolate a smaller number of key processes from the full climate system. They are typically formulated as systems of ordinary differential equations (ODEs) and allow for a comprehensive analysis within the framework of dynamical systems theory~\cite{kuznetsov1998elements, strogatz2001nonlinear, wiggins2003introduction}.

Conceptual box models have provided foundational insights into AMOC dynamics. Stommel's seminal two-box model~\cite{stommel1961thermohaline}, comprising a tropical and a subpolar North Atlantic box, formalized the salt--advection feedback: stronger overturning transports warmer, saltier water to high latitudes, increasing the North Atlantic surface density and further reinforcing the circulation. The model exhibits bistability between a thermally dominated state with northward overturning and a salinity-dominated state with reversed overturning. Alternatively, Welander~\cite{welander1982simple} formulated a two-box water-column model with vertical exchange represented by convective adjustment~\cite{ahmed2020deep, kurtze2010convective}. The model produces an on--off ``flip-flop'' cycle: under stable stratification, vertical exchange is weak and diffusive; once the column becomes statically unstable, convection mixes the boxes rapidly, reducing the density contrast and reestablishing stable stratification. Between convective events, background diffusion and surface fluxes slowly rebuild the vertical density difference~\cite{marshall1999open, rahmstorf2001simple}.

Several conceptual models extended those of Stommel and Welander by incorporating additional boxes and processes. For example, stochastic variants of Stommels model showed that bistability persists under random atmospheric forcing \cite{cessi1994simple,griffies1995linear} and that such noise can generate interdecadal oscillations \cite{rivin1997linear}. Interhemispheric multi-box frameworks established conditions for multiple equilibria under prescribed surface freshwater fluxes~\cite{huang1992multiple, rooth1982hydrology, welander1986thermohaline}, and related work coupled the ocean to simple energy balance models (EBM; see Figure~\ref{fig:climate_heirarcy})~\cite{birchfield1994century, birchfield1990coupled, lohmann1996stability}. Additionally, early conceptual studies advanced the ``salt--oscillator'' hypothesis~\cite{broecker1990salt, sakai1999dynamical}, in which salinity advection, freshwater forcing, and episodic convection interact to produce so-called Dansgaard--Oeschger variability~\cite{dokken2013dansgaard, li2019coupled}.

By the early--mid 2000s, dynamical systems theory was widely applied to low-order AMOC box models, and coupled salt-advection and convective-adjustment feedbacks had become commonplace, often producing millennial-scale relaxation oscillations. In particular, Colin de Verdi{\`e}re and co-authors~\cite{colin20002, colin2010instability, dijkstra2008dynamical}, demonstrated that slow salinity build-up followed by rapid convective ``flushes'' generates a sawtooth overturning cycle; their limit cycle is created by a Hopf bifurcation and terminates on a homoclinic bifurcation \cite{colin2007simple, de2006bifurcation}. Box models built for comparability with GCMs simulations -- most notably the Cimatoribus--Castellana model and the Stommel-Welander hybrid of Zhang and Marshall -- have been used to interpret features of coupled GCM simulations~\cite{castellana2019transition, cimatoribus2014meridional, soons2024optimal, van2024role, zhang2002mechanisms}. A complementary `top-down' route starts from a GCM, then constructs and calibrates a multi-basin box model to reproduce its overturning regimes~\cite{alkhayuon2019basin, smith2008description, wood2019observable}.

More recent extensions of conceptual AMOC box models introduce delay--differential equations (DDEs) to represent finite-time feedbacks loops, such as those with explicit advective travel times and convective adjustment memory~\cite{ManciniRuschelDijkstraKrauskopf2025TwoDelayAMOC}, thereby describing the AMOC dynamics with an infinite-dimensional phase space~\cite{erneux2009applied, smith2011introduction}. This class of equations preserves the low-order structure of box models while capturing delayed feedbacks in the overturning circulation, thereby reproducing complex dynamical phenomena more commonly seen in higher-complexity models~\cite{quinn2018delayed, wei2022simple}. Another approach examines rate-induced tipping (R-tipping) in the AMOC, whereby sufficiently rapid changes in time-varying external forcing result in transitions between different attractors without crossing a classical bifurcation (B-tipping)~\cite{ashwin2012tipping, feudel2018multistability, wieczorek2023rate}. For example, ice-sheet melt can trigger a tipping cascade \cite{klose2024rate}, resulting in a weakened AMOC. Basin stability and tipping thresholds~\cite{wieczorek2023rate} in a multi-ocean box model have been compared with those in the FAMOUS GCM~\cite{alkhayuon2019basin, chapman2024tipping}. Several studies also investigate AMOC recovery under time-dependent freshwater forcing near fold points, showing recovery if the forcing is reversed sufficiently rapidly \cite{budd2022dynamic, ritchie2021overshooting, sinet2024amoc, soons2025physics}.

The literature on AMOC box models is extensive, yet comprehensive bifurcation analyses are rare, and many studies rely on dimension-reducing assumptions that mask the full dynamical structure of the AMOC. To address this, we develop and analyze a six dimensional conceptual box model (later reduced to four dimensions) that combines and extends the classical frameworks of Stommel and Welander, largely following the setup of Zhang~\cite{zhang2002mechanisms}. We partition the ocean into three interconnected boxes: a tropical surface box, a subpolar surface box, and a deep-water box. Temperature and salinity in the surface boxes evolve through advective exchange between them, while both surface boxes interact with the deep-water box through convective-adjustment. The full model has six state variables (temperature and salinity in each box), which we reduce to four with timescale-motivated assumptions. This configuration retains the essential physics of advective and convective feedbacks while remaining tractable for a dynamical analysis.

We adopt a dynamical systems approach and carry out a comprehensive bifurcation analysis of our model. We examine its equilibria, representing  either a northward and southward steady state circulation, which may coexist in multistable parameter regimes. In addition, we identify periodic and chaotic solutions that generalise the \emph{Welander oscillations} recently documented in the classical two-dimensional setting of the Welander model~\cite{bailie2025detailed}. In our model, these oscillations switch between effectively `off' states, in which mixing is weak and diffusive, and strong predominantly convecting states with vigorous mixing. These transitions occur on a multi-decadal time scale that is slaved to a much slower, millennial-scale modulation of the AMOC strength. Notably, increasing the prescribed freshwater influx into the North Atlantic raises the frequency of convective shutdown events and weakens the overturning.

This paper is organized as follows. Section~\ref{section:model_formulation} formulates the six-dimensional three-box model for temperature and salinity and introduces two empirical assumptions that allow us to reduce this model to a system of four nondimensional equations. Section~\ref{section:bifurcation_structure} presents the equilibrium bifurcation structure, starting with a one-parameter bifurcation analysis and then extending it to a two-parameter setting. In this way, we identify the saddle--node and Hopf bifurcations that delimit several regions of multistability. This provides the foundation for Section~\ref{section:twopar_periodic_orbit}, where we analyze Welander oscillations in one- and two-parameter settings and characterize the associated pattern of convective shutdown events. Finally, Section~\ref{section:chaotic_pulse_welander} examines chaotic Welander oscillations, characterized by irregular convective shutdown events along chaotic trajectories; we also present here a partial analysis based on the maximal Lyapunov exponent that identifies the regions with chaotic dynamics. Taken together, these sections present a comprehensive dynamical systems analysis of our model, showing how equilibria, Welander oscillations, and chaotic shutdown dynamics are organized in a parameter plane of virtual salinity influx and a density threshold.

\section{Model formulation and background}
\label{section:model_formulation}
\begin{figure}[t!]
  \centering
  \includegraphics{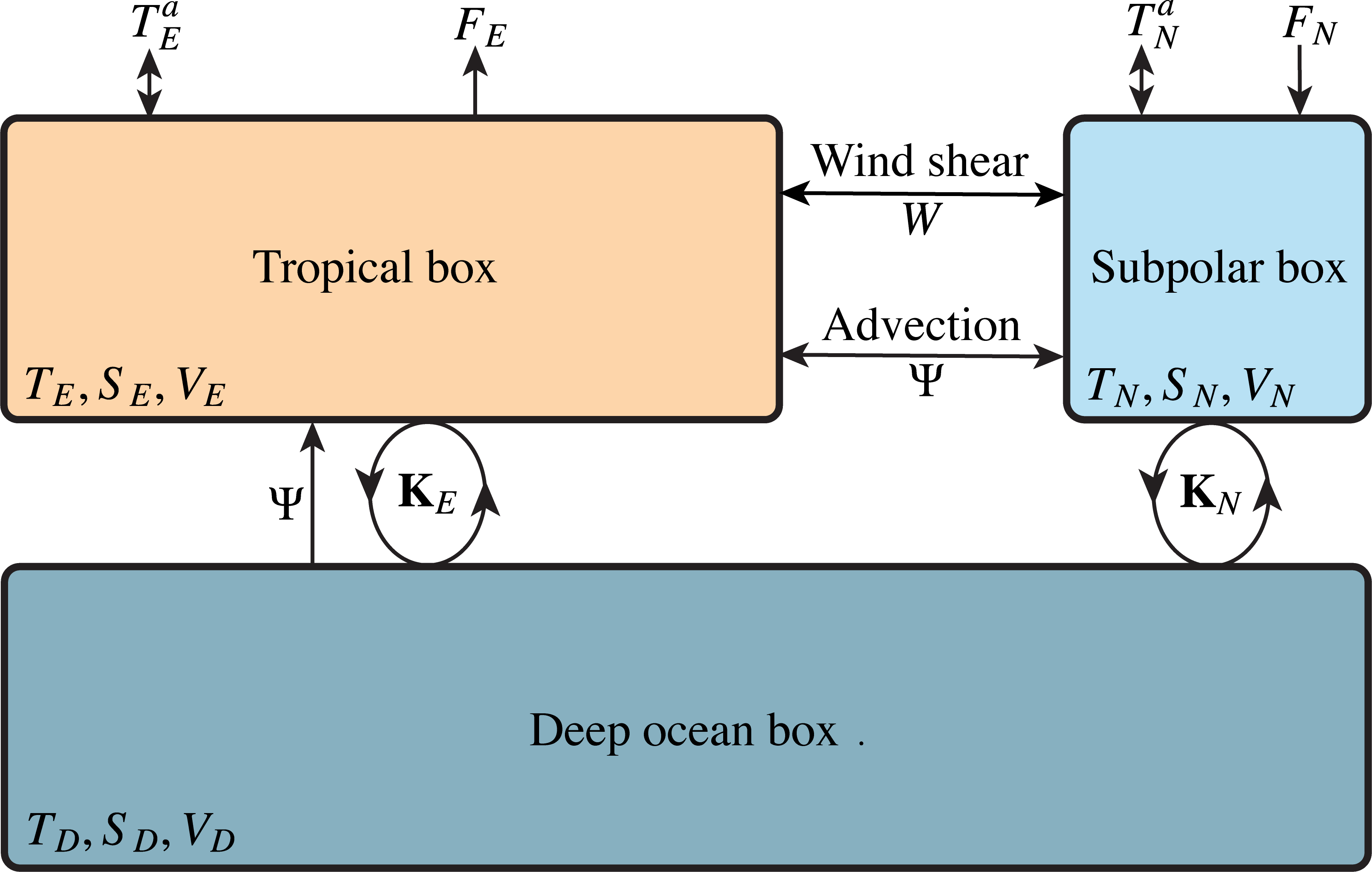}
  \caption{Model schematic for system~\eqref{eq:system_full}, with arrows indicating flux directions and loops indicating mixing, as governed by the convective exchange functions~$\mathbf{K}_E$ and~$\mathbf{K}_N$. The model consists of three interconnected basins: the tropical surface box $E$ (orange), the subpolar surface box $N$ (light blue), and the deep-water box $D$ (dark blue). Each box is labeled by its volume $V_i$, temperature $T_i$, and salinity $S_i$. In the tropical box, temperature relaxes toward the atmospheric value $T_E^a$ and salinity is forced by the virtual salt flux $F_E$. In the subpolar box, temperature relaxes toward $T_N^a$ and salinity is forced by the virtual salt flux $F_N$. Meridional exchange between the tropical and subpolar boxes is driven by a density-dependent overturning function~$\Psi$ and a wind-driven advection~$W$.}
  \label{fig:boxModelSchem}
\end{figure}

We construct a three-box configuration consisting of a tropical surface box, a subpolar surface box, and a deep box. Shown in Figure~\ref{fig:boxModelSchem}, each box is characterized by its temperature $T_i$, salinity $S_i$, and volume $V_i$, where $i \in \{E,N,D\}$ indexes the tropical (E), subpolar North Atlantic (N), and deep-water (D) boxes; arrows indicate the inflow and outflow associated with each box. The governing equations follow from heat and salt budgets under the assumption of fixed volumes $V_i$. All advective terms are written by using an upwind differencing scheme~\cite{calmanti2006north, rivin1998linear, zhang2002mechanisms}, and advection may be either northward or southward.

We adopt a linear equation of state~\cite{born2014two} for the box densities
\begin{equation}
  \rho_i= -\alpha_T\,(T_i-T_0)+\alpha_S\,(S_i-S_0), \qquad i\in\{E,N,D\},
\end{equation}
where $\alpha_T$ and $\alpha_S$ are constant thermal expansion and haline contraction coefficients, and $T_0$ and $S_0$ are reference values for the temperature and salinity~\cite{dijkstra2013nonlinear}.

Two advective processes are represented in Figure~\ref{fig:boxModelSchem}. The first is a wind--driven gyre transport, modeled as a linear meridional exchange of temperature and salinity with strength $W$~\cite{alkhayuon2019basin,longworth2005ocean}. The second is density-driven meridional exchange between the tropics and the subpolar North Atlantic. We parameterize the overturning with a Stommel--type closure~\cite{dijkstra2008dynamical, stommel1961thermohaline}
\begin{equation}
  \Psi = \sigma\,(\rho_N - \rho_E),
\end{equation}
where $\sigma>0$ is a hydraulic constant. Here $\Psi>0$ corresponds to northward surface flow with a southward deep return, while $\Psi<0$ indicates southward surface flow with northward deep return~\cite{stommel1961thermohaline}. Upwelling from the deep box to the surface boxes scales with the overturning strength $|\Psi|$ and is routed to the tropical or subpolar box according to the sign of $\Psi$. To ensure smoothness near $\Psi = 0$, we regularize the absolute value $|\Psi|$ as
\begin{equation}
  \mathbf \Psi = \Psi\,\tanh\!\Big(\frac{\Psi}{\vartheta}\Big),
  \qquad 0<\vartheta\ll 1,
  \label{eq:psi_regularization}
\end{equation}
which is $C^{\infty}$-smooth for any $\vartheta > 0$ and converges uniformly to $|\Psi|$ as $\vartheta \to 0$.

We parameterize vertical exchange between each surface box and the deep water box by a Welander-type convective adjustment~\cite{bailie2025detailed, bailie2024bifurcation, welander1982simple}. For $i \in \{N,E\}$, we introduce the convective exchange function
\begin{equation}
  \mathbf K_{i} = k_d + \frac{1}{2}\,(k_i - k_d)\left(1 + \tanh\!\left(\frac{\rho_i - \rho_D - \Delta\rho}{\varepsilon}\right)\right),
  \label{eq:conv_ex_subpolar}
\end{equation}
where $k_d>0$ is the background diffusive coupling and $k_{i}>k_d$ is the convective mixing strength in box $i$. The parameter $\varepsilon>0$ controls the sharpness of the transition between the diffusive and convective regimes. The interpretation is that background diffusion maintains a weak but continuous coupling, while convection initiates only when the density contrast $\rho_i - \rho_D$ between the surface and deep water is sufficiently larger than the density threshold $\Delta\rho$. We impose that the mixing strengths satisfy $0<k_d<k_E<k_N$, which is consistent with stronger subpolar convection and more stable tropical stratification~\cite{zhang2002mechanisms}.

In the surface boxes, exchange with an effectively infinite atmospheric reservoir supplies the remaining fluxes. Tropical temperature $T_E$ is relaxed to its prescribed atmospheric value $T^a_E$ at rate $\gamma_E$ over volume $V_E$, and salinity receives a virtual salt flux $F_E S_0$ representing net evaporation. In the subpolar North Atlantic, $T_N$ is restored to $T^a_N$ at rate $\gamma_N$ over $V_N$, and salinity is forced by a virtual salinity flux $F_N S_0$ representing net freshwater input from precipitation and Greenland Ice Sheet runoff~\cite{bamber2012recent}. We adopt the convention that $F_i > 0$ increases salinity (evaporation), and $F_i<0$ decreases salinity (freshwater input).

Combining the oceanic exchange and atmospheric forcing, the temperature and salt budgets in the tropical (E), subpolar (N), and deep (D) boxes evolve according to
\begin{equation}\label{eq:system_full}
  \begin{aligned}
    V_N \frac{dT_N}{dt} &= -\gamma_N V_N (T_N - T_N^a)
    + \frac{\Psi}{2}(T_E - T_D)
    + W(T_E - T_N)
    + \frac{|\Psi|}{2}(T_E + T_D - 2T_N)  \\
    &\quad
    - \mathbf{K}_N\,(T_N - T_D),                   \\[1ex]
    V_N \frac{dS_N}{dt} &= F_N S_0
    + \frac{\Psi}{2}(S_E - S_D)
    + \frac{|\Psi|}{2}(S_E + S_D - 2S_N)
    + W(S_E - S_N)  \\
    &\quad
    + \mathbf{K}_N\,(S_D - S_N),                   \\[1ex]
    V_E \frac{dT_E}{dt} &= -\gamma_E V_E (T_E - T_E^a)
    + \frac{\Psi}{2}(T_D - T_N)
    + \frac{|\Psi|}{2}(T_D + T_N - 2T_E)
    - W(T_E - T_N) \\
    &\quad
    - \mathbf{K}_E\,(T_E - T_D),                    \\[1ex]
    V_E \frac{dS_E}{dt} &= F_E S_0
    + \frac{\Psi}{2}(S_D - S_N)
    + \frac{|\Psi|}{2}(S_D + S_N - 2S_E)
    - W(S_E - S_N)  \\
    &\quad
    - \mathbf{K}_E\,(S_E - S_D),                    \\[1ex]
    V_D \frac{dT_D}{dt} &= \frac{\Psi}{2}(T_N - T_E)
    + \frac{|\Psi|}{2}(T_N + T_E - 2T_D)
    + \mathbf{K}_N\,(T_N - T_D)  \\
    &\quad
    + \mathbf{K}_E\,(T_E - T_D),                     \\[1ex]
    V_D \frac{dS_D}{dt} &= \frac{\Psi}{2}(S_N - S_E)
    + \frac{|\Psi|}{2}(S_N + S_E - 2S_D)
    + \mathbf{K}_N\,(S_N - S_D)  \\
    &\quad
    + \mathbf{K}_E\,(S_E - S_D).
  \end{aligned}
\end{equation}
This six-dimensional system fully describes the temperature and salinity dynamics of the three-box configuration in~Figure~\ref{fig:boxModelSchem}.

\subsection{Dimension reduction and non-dimensionalisation}
\label{section:model_reduction}

We now assume that: (i) the deep-ocean salinity $S_D$ evolves on a much longer timescale and is, therefore, held constant~\cite{alkhayuon2019basin}, and (ii) that the tropical temperature $T_E$ relaxes rapidly to its prescribed atmospheric value $T_E^a$ on a very fast timescale, so that we set $T_E = T_E^a$. The resulting four-dimensional reduced model is given by
\begin{equation}\label{eq:system_reduced}
  \begin{aligned}
    V_N \frac{dT_N}{dt} &= -\gamma_N V_N (T_N - T_N^a)
    + \frac{\Psi}{2}(T_E^a - T_D)
    + W(T_E^a - T_N)
    + \frac{|\Psi|}{2}(T_E^a + T_D - 2T_N)  \\
    &\quad
    - \mathbf{K}_N\,(T_N - T_D),                   \\[1ex]
    V_N \frac{dS_N}{dt} &= F_N S_0
    + \frac{\Psi}{2}(S_E - {S}_D)
    + \frac{|\Psi|}{2}(S_E + {S}_D - 2S_N)
    + W(S_E - S_N)  \\
    &\quad
    + \mathbf{K}_N\,({S}_D - S_N),                   \\[1ex]
    V_E \frac{dS_E}{dt} &= -F_N S_0
    + \frac{\Psi}{2}({S}_D - S_N)
    + \frac{|\Psi|}{2}({S}_D + S_N - 2S_E) - W(S_E - S_N)  \\
    &\quad
    - \mathbf{K}_E\,(S_E - {S}_D),                    \\[1ex]
    V_D \frac{dT_D}{dt} &= \frac{\Psi}{2}(T_N - T_E^a)
    + \frac{|\Psi|}{2}(T_N + T_E^a - 2T_D)
    + \mathbf{K}_N\,(T_N - T_D)  \\
    &\quad
    + \mathbf{K}_E\,(T_E^a - T_D),                     \\[1ex]
  \end{aligned}
\end{equation}
where we impose a zero net freshwater influx, such that $F_N = -F_E$.

\begin{table}[b!]
  \centering
  \caption{Physical parameters and their values used in the three-box ocean model. Salinity is dimensionless. Flow parameters are given in Sv, where $1~\mathrm{Sv}=10^6~\mathrm{m}^3/\mathrm{s}$.}
  \label{table1:physicalParm}
  \begin{tabular}{@{}p{0.60\linewidth} l r l@{}}
    \toprule
    \textbf{Parameter} & \textbf{Symbol} & \textbf{Value} & \textbf{Units} \\
    \midrule

    \multicolumn{4}{@{}l}{\textit{Physical Constants}} \\
    Temperature expansion coefficient & $\alpha_T$ & $1.7 \times 10^{-4}$ & $^{\circ}$C$^{-1}$ \\
    Salinity contraction coefficient  & $\alpha_S$ & $0.8 \times 10^{-3}$ & \\
    North Atlantic atmospheric temperature & $T^a_N$ & $7$ & $^{\circ}$C \\
    Equatorial atmospheric temperature     & $T^a_E$ & $25$ & $^{\circ}$C \\
    Reference temperature                  & $T_0$   & $2.65$ & $^{\circ}$C \\
    Reference salinity                     & $S_0$   & $35$ & -- \\
    Deep water salinity                    & $S_D$   & $34.538$ & -- \\
    Relaxation coefficient                   & $\gamma_N$   & $1$ & $yr^{-1}$ \\
    \midrule

    \multicolumn{4}{@{}l}{\textit{Box Volumes}} \\
    North Atlantic box volume & $V_N$ & $7.2106 \times 10^{15}$ & m$^3$ \\
    Equatorial box volume     & $V_E$ & $6.3515 \times 10^{16}$ & m$^3$ \\
    Deep water box volume     & $V_D$ & $1.2979 \times 10^{17}$ & m$^3$ \\
    \midrule

    \multicolumn{4}{@{}l}{\textit{Flow Parameters}} \\
    Hydraulic constant             & $\sigma$   & $2.1 \times 10^{4}$ & Sv \\
    Wind advection strength        & $W$        & $5.456$             & Sv \\
    Subpolar convective strength   & $\kappa_N$ & $8.0$               & Sv \\
    Equatorial convective strength & $\kappa_E$ & $4.0$               & Sv \\
    Global diffusivity             & $\kappa_d$ & $0.01$              & Sv \\
    \bottomrule
  \end{tabular}
\end{table}
To non-dimensionalize system~\eqref{eq:system_reduced}, we introduce the rescaled parameters
\begin{align*}
  \tilde t = &\frac{\sigma}{V_N}t \qquad \mu = \frac{F_N S_0\,\alpha_S}{\sigma\,\alpha_T\,(T_N^a - T_0)},
  \qquad
  \eta   = \frac{\Delta\rho}{\alpha_T\,(T_N^a - T_0)\,\rho_0}, \qquad \phi = \alpha_T (T^a_N - T_0),
  \\[0.5ex]
  \widetilde{W} &= \frac{W}{\sigma}, \qquad
  \widetilde{V}_i = \frac{V_i}{V_N}, \qquad
  \kappa_{c/d}^i = \frac{k_{c/d}^i}{\sigma}, \qquad
  \delta_N = \frac{\gamma_N}{\sigma}.
\end{align*}
Here $\tilde t$ sets the characteristic time scale through the hydraulic constant~$\sigma$ and the North Atlantic volume~$V_N$, the parameter $\mu$ measures the strength of the virtual freshwater forcing relative to the advective salt transport into the North Atlantic, and $\eta$ is the corresponding dimensionless density threshold for the onset of convection. The resulting non-dimensionalized overturning strength (with some slight abuse of notation) is given by
\begin{equation}
  \Psi = \phi \bigl[ (y_N - x_N) - (y_E - x_E^a) \bigl],
\end{equation}
and the non-dimensionalized convective exchange functions are
\begin{align}
  \mathbf{K}_N &= \kappa_d
  + \tfrac{1}{2}\big(\kappa_N - \kappa_d\big)
  \left(1 + \tanh\!\left(\frac{(y_N - x_N) - (y_D - x_D) - \eta}{\varepsilon}\right)\right),
  \label{eq:conv_ex_subpolar}
  \\[1ex]
  \mathbf{K}_E &= \kappa_d
  + \tfrac{1}{2}\big(\kappa_E - \kappa_d\big)
  \left(1 + \tanh\!\left(\frac{(y_E - x_E^a) - (y_D - x_D) - \eta}{\varepsilon}\right)\right).
  \label{eq:conv_ex_tropical}
\end{align}

Substituting these scalings into system~\eqref{eq:system_reduced} leads to the non-dimensional reduced model
\begin{equation}
  \begin{aligned}
    \dot x_N &= -\delta_N(x_N-1)
    + \frac{\Psi}{2}(x_E^a - x_D)
    + \frac{|\Psi|}{2}(x_E^a + x_D - 2x_N)
    + \widetilde{W}(x_E^a - x_N) \\
    &\quad
    - \mathbf{K}_N(x_N - x_D),
    \\[1ex]
    \dot y_N &= \mu
    + \frac{\Psi}{2}(y_E - y_D)
    + \frac{|\Psi|}{2}(y_E + y_D - 2y_N)
    + \widetilde{W}(y_E - y_N)
    \\
    &\quad
    - \mathbf{K}_N(y_N - y_D),
    \\[1ex]
    \widetilde V_E\dot y_E &= \!
    -\mu
    + \frac{\Psi}{2}(y_D - y_N)
    + \frac{|\Psi|}{2}(y_D + y_N - 2y_E)
    - \widetilde{W}(y_E - y_N)  \\
    &\quad
    - \mathbf{K}_E(y_E - y_D),
    \\[1ex]
    \widetilde V_D\dot x_D &= \!
    \frac{\Psi}{2}(x_N - x_E^a)
    + \frac{|\Psi|}{2}(x_N + x_E^a - 2x_D)
    + \mathbf{K}_N(x_N - x_D)  \\
    &\quad
    + \mathbf{K}_E(x_E^a - x_D).
  \end{aligned}
  \label{eq:non_dim_system}
\end{equation}
Here, the dot denotes differentiation with respect to the rescaled time~$\tilde t$. In our scaling, one unit of~$\tilde t$ corresponds to approximately $3.97$ days in dimensional time~$t$. For notational convenience, we drop the tilde and write simply $t$ for the rescaled time from now on. To aid in the physical interpretation of the model, Table~\ref{table1:physicalParm} lists the dimensional parameters assembled from previous box-model studies and observational data.

The four-dimensional scaled system~\eqref{eq:non_dim_system} is our object of study. We treat the North Atlantic virtual salinity $\mu$ and the density threshold $\eta$ as the primary bifurcation parameters. Following standard practice, we take $\eta<0$ so that convective adjustment is triggered before the onset of static instability in the water column~\cite{bailie2025detailed, bailie2024bifurcation, cessi1994simple, welander1982simple}. We focus on regimes with a net freshwater influx into the North Atlantic, that is, $\mu<0$. The bifurcation structure depends sensitively on the smoothing parameter $\varepsilon$ and shifts significantly as $\varepsilon \to 0$. Moreover, the parameters~$\mu$ and~$\eta$ are phenomenological, in the sense that they represent effective controls that aggregate unresolved climate processes driving advection and preconditioning convective mixing. For these two reasons, we allow a broader range of~$(\mu,\eta)$ than would be suggested by a direct physical interpretation. Throughout, we implement a smooth switching with $\varepsilon = 0.02$ in~\eqref{eq:conv_ex_subpolar} and~\eqref{eq:conv_ex_tropical}; this choice remains numerically tractable and lies in a parameter range where Welander oscillations are observed in the Welander model~\cite{bailie2025detailed}. For consistency, we also choose $\vartheta = 0.02$ as the smoothing parameter in~\eqref{eq:psi_regularization}.

\section{Bifurcation structure of equilibria}
\label{section:bifurcation_structure}
We now analyze the equilibrium structure of system~\eqref{eq:non_dim_system} as the virtual salinity flux $\mu$ and density threshold $\eta$ are varied. Our goal is to identify all possible steady overturning configurations and their stability. To this end, we first continue equilibria in $\mu$ for a selection of fixed values of~$\eta$ and then extend this analysis to the full two-parameter $(\mu,\eta)$-plane. Curves $\mathrm{SN}$ of saddle--node bifurcations and $\mathrm{H}$ of Hopf bifurcations partition the parameter plane into open regions with different numbers of stable equilibria, corresponding to distinct northward and southward overturning AMOC states.

\subsection{One-parameter bifurcation analysis}
\label{section:onepar_steady_state}
\begin{figure}[t!]
  \centering
  \includegraphics{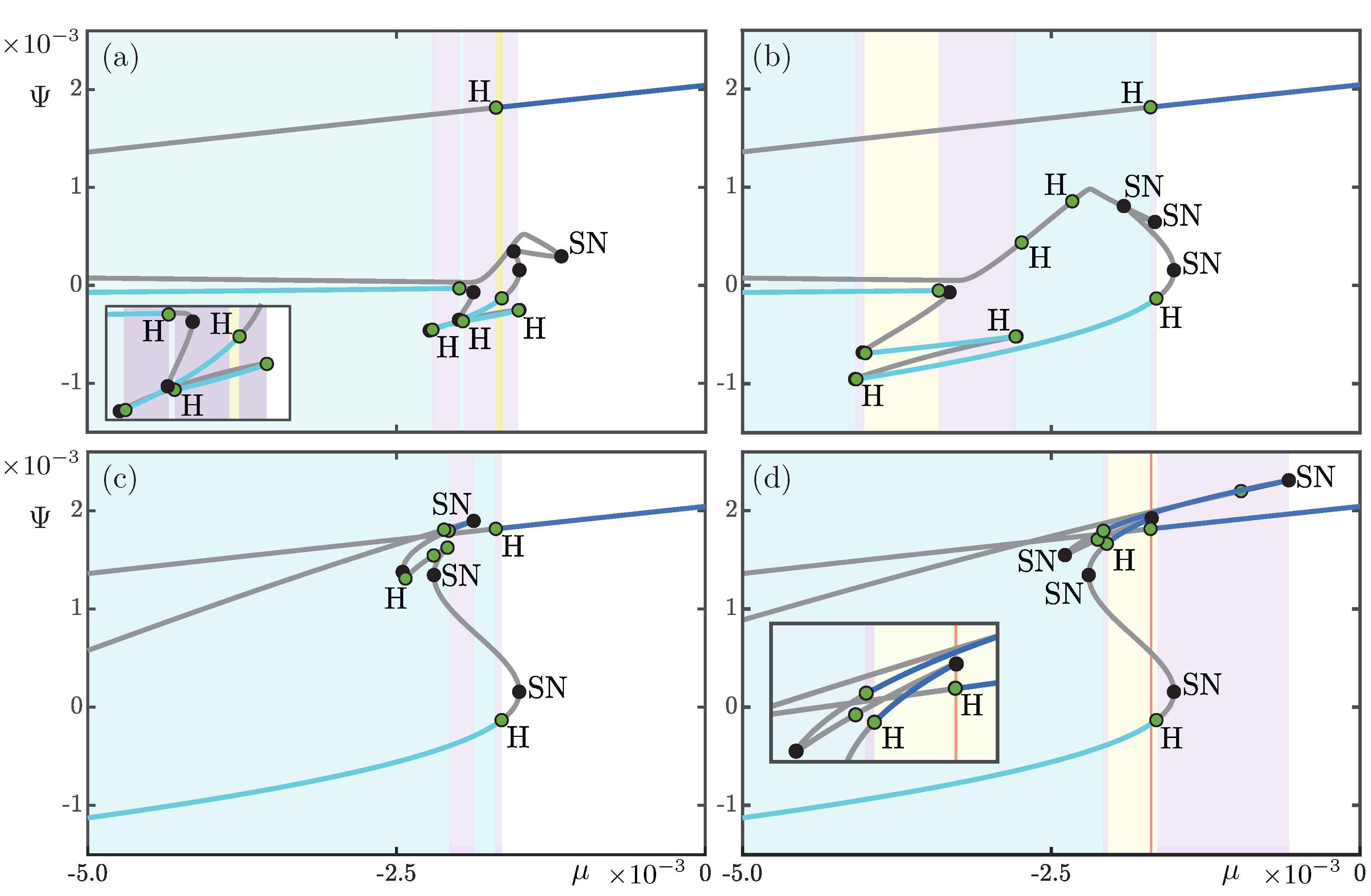}
  \caption{One-parameter bifurcation diagrams in the freshwater forcing~$\mu$, showing branches of equilibria represented by the overturning function~$\Psi$. Points $\mathrm{SN}$ of saddle--node and $\mathrm{H}$ of Hopf bifurcations are marked by black and green circles. The branches are colored to indicate their stability: stable northward (blue), stable southward (cyan), and unstable (black) equilibria. Background shading denotes $\mu$-intervals with one northward equilibrium (unshaded), one southward equilibrium (cyan), and two (purple), three (yellow), or four (red) coexisting stable equilibria. Panels~(a)--(d) are for $\eta=-0.5$, $\eta=-1.5$, $\eta=-3.99$, and $\eta=-5.0$, respectively.} \label{fig:EqCurves}
\end{figure}

Figure~\ref{fig:EqCurves} shows branches of equilibria of system~\eqref{eq:non_dim_system}, for four successively smaller values of the density threshold~$\eta$, continued in the virtual salinity~$\mu$ and plotted against the overturning function~$\Psi$ from~\eqref{eq:non_dim_system}. Stable equilibria with northward overturning are shown in blue, those with southward overturning in cyan, and unstable equilibria in black. The background shading  of Figure~\ref{fig:EqCurves} indicates the number of coexisting \emph{stable} equilibria: unshaded for a single northward state (positive $\Psi$), cyan for a single southward state (negative $\Psi$), and purple, yellow, and red for two, three, and four coexisting stable states, respectively. Note that, the background shading distinguishes between northward and southward only in monostable regimes, while in multistable regimes it encodes solely the number of coexisting stable equilibria. In each panel of Figure~\ref{fig:EqCurves}, we begin with $\mu$ near zero, where a single stable northward-overturning branch exists; we then trace how the representative equilibrium structure changes as $\mu$ is decreased.

Figure~\ref{fig:EqCurves}(a), for $\eta = -0.5$, shows that up to three stable equilibrium branches can coexist. For $\mu$ near zero, the unshaded interval contains a single stable branch with comparatively large northward overturning $\Psi > 0$. As $\mu$ decreases, a southward-overturning branch gains stability at a Hopf bifurcation~$\mathrm{H}$, marking the onset of an interval of bistability (purple) with the northward branch. A narrow yellow interval of tristability emerges at a second Hopf bifurcation~$\mathrm{H}$, where a second southward-overturning branch becomes stable. The northward branch then loses stability at a Hopf bifurcation $\mathrm{H}$, after which the purple interval corresponds to a regime with two coexisting stable southward equilibria. For more negative~$\mu$, one of these equilibria losses stability at a Hopf bifurcation $\mathrm{H}$, and system~\eqref{eq:non_dim_system} enters the cyan regime, where only a single stable southward equilibrium exists; note that $\Psi < 0$ is almost zero in absolute value.

Figure~\ref{fig:EqCurves}(b), for $\eta = -1.05$, shows a similar onset of multistability. As $\mu$ decreases, a southward-overturning branch gains stability at a Hopf bifurcation~$\mathrm{H}$, creating a narrow purple interval of bistability with the northward branch. At a second Hopf bifurcation point~$\mathrm{H}$, the northward branch loses stability, and system~\eqref{eq:non_dim_system} enters a broad cyan regime with a single stable southward equilibrium. For more negative~$\mu$, the southward equilibrium family develops alternating stable and unstable segments. Their stability changes at successive Hopf bifurcation points~$\mathrm{H}$, and saddle--node bifurcations~$\mathrm{SN}$ mark folds of the branch. The resulting sequence of stable intervals terminates, leaving a single stable southward-overturning equilibrium in the cyan regime of $\mu$.

Figure~\ref{fig:EqCurves}(c), for $\eta=-3.99$, shows a different sequence with only monostable or bistable regimes. As $\mu$ decreases from the unshaded regime, a southward-overturning branch gains stability at a Hopf bifurcation~$\mathrm{H}$, producing a narrow purple interval of bistability with the northward branch. At a subsequent Hopf point~$\mathrm{H}$, the northward branch loses stability and system~\eqref{eq:non_dim_system} enters a cyan regime with a single stable southward equilibrium. For more negative~$\mu$, a new northward branch turns at a saddle--node bifurcation~$\mathrm{SN}$ and becomes stable, creating a second purple interval in which it coexists with a stable southward equilibrium. This northward branch then loses stability again, after which system~\eqref{eq:non_dim_system} enters the leftmost cyan regime with only a single stable southward equilibrium.

Figure~\ref{fig:EqCurves}(d), for $\eta = 5.0$, shows a more intricate sequence of multistable regimes. As $\mu$ decreases from the unshaded regime, system~\eqref{eq:non_dim_system} first enters a broad purple interval in which two northward equilibria coexist. It then passes through a narrow yellow interval and an even narrower red interval, bounded by a saddle--node bifurcation~$\mathrm{SN}$ and a Hopf bifurcation~$\mathrm{H}$. For more negative~$\mu$, the sequence of multistable regimes reverses: system~\eqref{eq:non_dim_system} passes back through yellow and then purple intervals before entering a cyan regime with a single stable southward equilibrium, which persists for all more negative~$\mu$.

Taken together, Figure~\ref{fig:EqCurves} shows that, in each panel, a single stable northward-overturning branch with comparatively large $\Psi > 0$ exists for $\mu$ near zero, whereas for sufficiently negative~$\mu$ only a stable southward equilibrium with $\Psi<0$ persists. The four cases differ in the locations and widths of the bi- and multistable intervals that arise as Hopf bifurcations~$\mathrm{H}$ and saddle--node bifurcations~$\mathrm{SN}$ are crossed. In this way, Figure~\ref{fig:EqCurves} illustrates how increasing freshwater forcing (decreasing~$\mu$) both weakens the AMOC and promotes parameter regimes with multiple coexisting steady overturning states.

\subsection{Equilibrium structure in the $(\mu,\eta)$-plane}
\label{section:twopar_steady_state}

\begin{figure}[t!]
  \centering
  \includegraphics{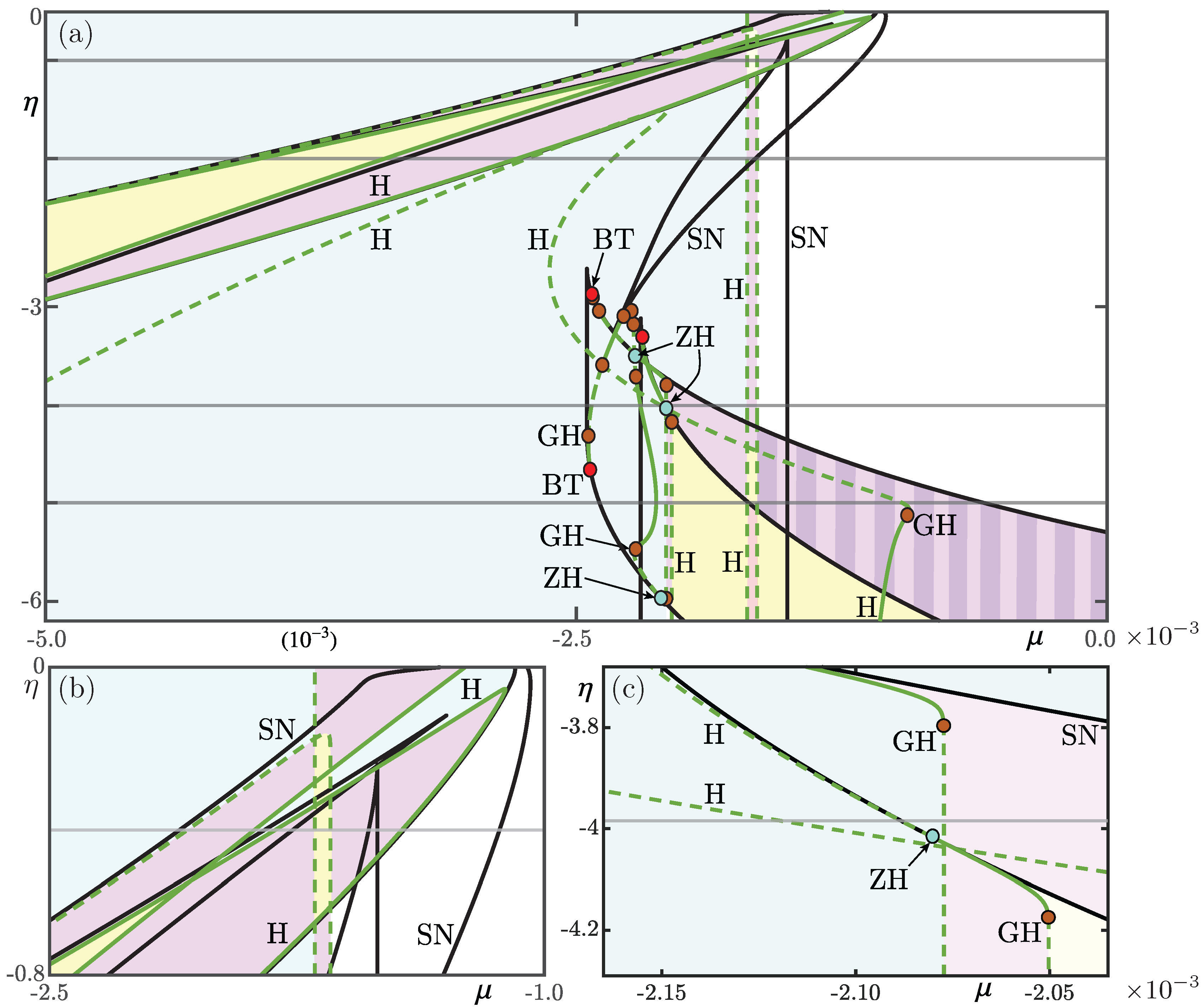}
  \caption{\label{fig:EqTwoPar}
  Two-parameter bifurcation diagram of equilibria in the $(\mu,\eta)$-plane of system~\eqref{eq:non_dim_system}. Panel~(a) shows a large portion of the parameter plane, with two enlargements in panels~(b) and~(c). Also shown are curves $\mathrm{SN}$ of saddle--node bifurcation (grey) and curves $\mathrm{H}$ of Hopf bifurcation (red; solid for supercritical and dashed for subcritical). Codimension-two points along these curves are indicated: Bogdanov--Takens $\mathrm{BT}$ (red), generalized Hopf $\mathrm{GH}$ (brown), and zero--Hopf $\mathrm{ZH}$ (cyan). The curves $\mathrm{SN}$ and $\mathrm{H}$ delineate regions with different numbers and types of stable equilibria: one northward equilibrium (unshaded); one southward equilibrium (cyan); two equilibria (purple); three equilibria (yellow); and four equilibria (red). The grey lines correspond to $\eta$-values used in Figure~\ref{fig:EqCurves}.}
\end{figure}
Figure~\ref{fig:EqTwoPar} presents the two-parameter equilibrium bifurcation structure of system~\eqref{eq:non_dim_system} in the $(\mu,\eta)$-plane. Panel~(a) illustrates how curves of Hopf~$\mathrm{H}$ and saddle--node~$\mathrm{SN}$ bifurcations partition the relevant part of the $(\mu,\eta)$-plane into open regions with different configurations of stable equilibria. Here, the shading follows the convention of Figure~\ref{fig:EqCurves}, and textured overlays are added to distinguish adjacent regions that have the same number of stable equilibria but differ in overturning type. Horizontal gray lines in Figure~\ref{fig:EqTwoPar} mark the four $\eta$-values used in the one-parameter continuations of Figure~\ref{fig:EqCurves}. Panels~(b) and~(c) enlarge neighborhoods where $\mathrm{H}$ and $\mathrm{SN}$ bound particularly narrow regions. Along the bifurcation curves, we also identify codimension--two points, including generalized Hopf (Bautin)~$\mathrm{GH}$, Bogdanov--Takens~$\mathrm{BT}$, and zero--Hopf~$\mathrm{ZH}$ bifurcations. In particular, the Hopf curves~$\mathrm{H}$ emerge from points~$\mathrm{BT}$ and change criticality at points~$\mathrm{GH}$ and~$\mathrm{ZH}$.

Consistent with the one-parameter diagrams in Figure~\ref{fig:EqCurves}, a broad unshaded region occupies values of~$\mu$ closer to zero in Figure~\ref{fig:EqTwoPar}, where a single northward equilibrium is stable. As~$\eta$ decreases from this region, one encounters a vertically striped purple region that indicates bistability between two northward-overturning equilibria. Adjacent lie three disjoint yellow regions of tristability, where the two northward equilibria coexist with a southward-overturning equilibrium that becomes stable upon crossing the relevant $\mathrm{SN}$ and $\mathrm{H}$ curves. Separating the yellow regions are two nearly vertical Hopf bifurcation curves~$\mathrm{H}$ and a saddle--node bifurcation curve~$\mathrm{SN}$, which together carve out a thin strip in which an additional northward equilibrium is stable; this strip corresponds to the narrow red $\mu$-interval in Figure~\ref{fig:EqCurves}(d). For more negative values of~$\mu$, two cyan regions exist that are separated by a long belt of alternating purple and yellow regions, where the southward branch coexists with one or two northward equilibria, respectively.

The narrow regions in panels~(b) and~(c) are of particular interest, since here the curves $\mathrm{SN}$ and $\mathrm{H}$ accumulate and bound intricate configurations of stable equilibria. Panel~(b) magnifies a small region for relatively small $\eta$, where $\mathrm{SN}$ and $\mathrm{H}$ curves cluster together to create small pockets of bistability and tristability. Panel~(c) focuses on an intermediate $\eta$-range in a neighborhood of a codimension--two zero--Hopf bifurcation~$\mathrm{ZH}$, where $\mathrm{SN}$ and $\mathrm{H}$ meet tangentially; this bifurcation acts as an organizing center for the surrounding equilibrium structure. The resulting arrangement closely mirrors the one-parameter picture in Figure~\ref{fig:EqCurves}(c): two disjoint monostable cyan regions are separated by a purple region of bistability, and a smaller yellow region of tristability appears in the lower-right corner.

Overall, decreasing the virtual salinity~$\mu$ drives system~\eqref{eq:non_dim_system} toward southward-overturning equilibria and progressively introduces bistable and tristable regimes. This behaviour is consistent with the salt--advection feedback in classical box models~\cite{lucarini2005thermohaline, rahmstorf1996freshwater} and with hosing experiments in OGCMs~\cite{manabe1995simulation,stouffer2006investigating,vellinga2002global}, which likewise predict weakened and potentially multiple AMOC steady states under enhanced freshwater forcing.

\section{Bifurcation analysis of Welander oscillations}
\label{section:twopar_periodic_orbit}
Having established how Hopf~$\mathrm{H}$ and saddle--node~$\mathrm{SN}$ bifurcations organize the equilibrium structure in the $(\mu,\eta)$-plane, we now turn to periodic solutions of system~\eqref{eq:non_dim_system}. Specifically, we introduce \emph{Welander oscillations}: self-sustained oscillations generated by the convective-adjustment scheme~$\mathbf{K}_N$, in which North Atlantic mixing alternates between convective and diffusive phases. The term originates from the classical two-dimensional Welander model, where the relevant periodic orbit has period one and features a single convective shutdown event~\cite{bailie2025detailed}. Here we consider the higher-dimensional setting, in which periodic solutions may exhibit multiple shutdown events and may also give rise to irregular, chaotic shutdown dynamics. Our goal is to analyze these oscillations and to identify the bifurcations that delimit their existence; we do not aim to reproduce a present-day climate state. Accordingly, we restrict attention to the smaller parameter range shown in Figure~\ref{fig:EqTwoPar}(c), focusing on a neighborhood of the codimension--two zero--Hopf point~$\mathrm{ZH}$, since such points act as organizing centers for families of periodic orbits and global bifurcations in parameter space~\cite{krauskopf2003saddle, krauskopf2004planar, zimmermann2001global}.

We carry out a combined one- and two-parameter bifurcation analysis: we first continue families of periodic solutions in the virtual salinity~$\mu$ for four successively decreasing values of the density threshold~$\eta$, and then extend this analysis to the $(\mu,\eta)$-plane, where we compute the loci of secondary saddle--node~$\mathrm{S}$ and period--doubling~$\mathrm{PD}$ bifurcations. These curves bound open parameter regions in which system~\eqref{eq:non_dim_system} exhibits stable Welander oscillations. For each such region, we present three complementary phase-space representations that highlight the pattern of convective shutdown events over one period; in particular, we find that the number of shutdown events per period increases as $\mu$ decreases.

\subsection{Welander oscillations}
To identify Welander oscillations and to quantify convective shutdown events along a periodic solution~$\Gamma$, we follow the geometric classification of~\cite{bailie2025detailed}. For $\varepsilon>0$, the convective-adjustment function~$\mathbf{K}_N$ is smooth and distinguishes three mixing regimes: a convective regime, a diffusive regime, and an intermediate regime created by regularization of the discontinuity at $\varepsilon=0$. We view the graph of~$\mathbf{K}_N$ as a surface over the $(\rho_B,\rho_N)$-plane, where $\rho_B$ denotes the deep-water density and $\rho_N$ the subpolar North Atlantic density, and define $L^{+}$ and $L^{-}$ to be the loci of maximal and minimal curvature of this surface. The intersections of~$\Gamma$ with $L^{+}\cup L^{-}$ then partition~$\Gamma$ into three segments, which we label as convective, diffusive, or intermediate according to the values of~$\mathbf{K}_N$ along each segment.

\begin{figure}[t!]
  \centering
  \includegraphics{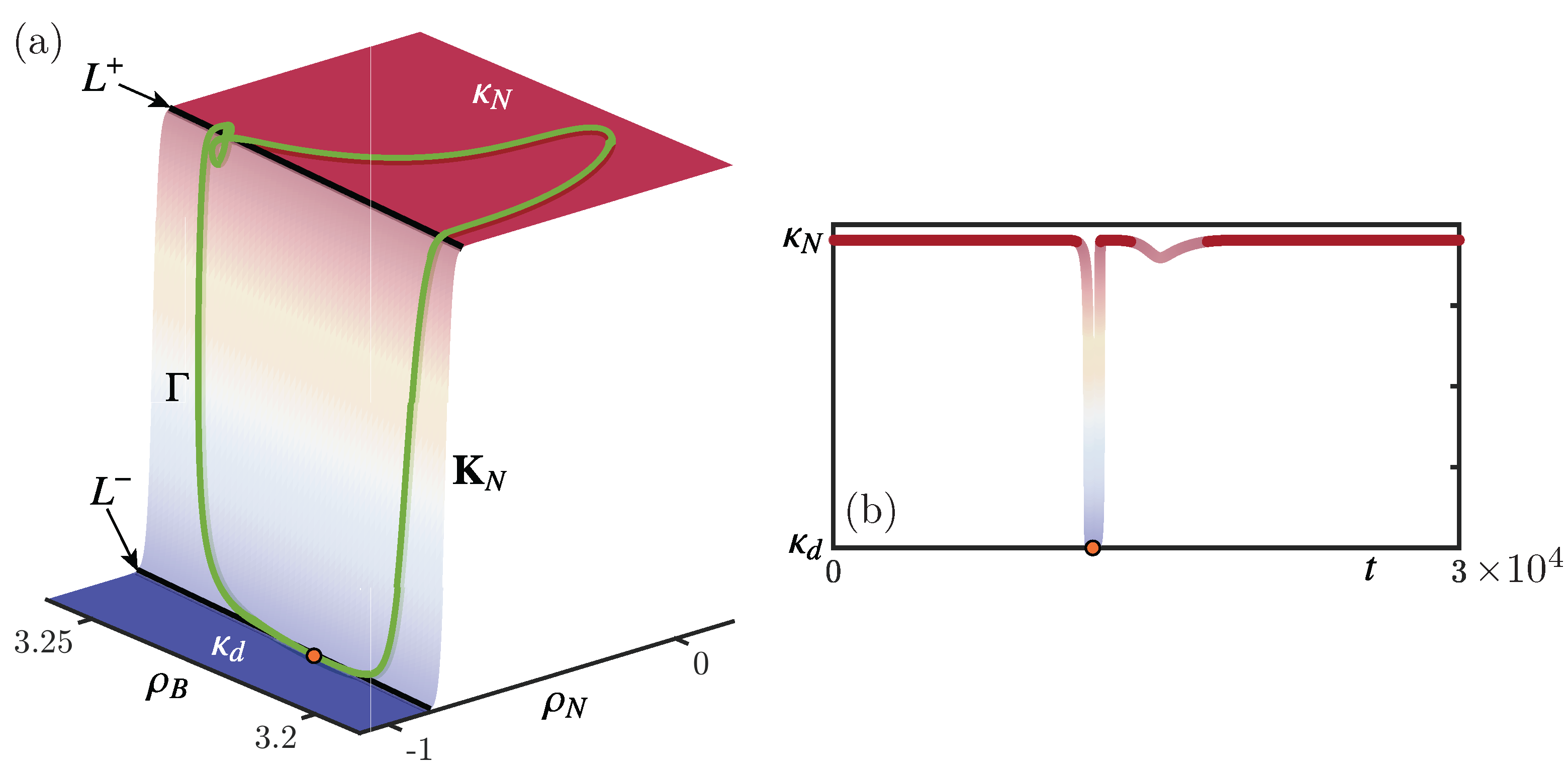}
  \caption{\label{fig:welander_occ_intro}
    Geometric classification of mixing regimes along a periodic orbit $\Gamma$ representing a Welander oscillation in system~\eqref{eq:non_dim_system} for $(\mu,\eta)=(-2.11\times10^{-3},-3.99)$. Panel~(a) shows $\Gamma$ in projection onto the surface of the convective exchange function~$\mathbf{K}_N$ over the $(\rho_B,\rho_N)$-plane, which is shaded by mixing regime: convective (red), diffusive (blue), and a narrow intermediate regime (semi-transparent). The boundaries between these regimes are the curvature loci $L^{+}$ and $L^{-}$ of maximal and minimal curvature. Panel~(b) shows the corresponding time series $\mathbf{K}_N(t)$ along~$\Gamma$; orange dots mark the minimum value of~$\mathbf{K}_N$ during each convective shutdown event.
  }
\end{figure}

Figure~\ref{fig:welander_occ_intro} illustrates this phase decomposition in two complementary representations. Panel~(a) shows the periodic solution~$\Gamma$ displayed on the surface of~$\mathbf{K}_N$ with the curvature loci $L^{+}$ and $L^{-}$; crossings of~$\Gamma$ through either locus mark transitions between mixing regimes. The upper sheet of the surface, where $\mathbf{K}_N\approx \kappa_N$, corresponds to the convective regime, whereas the lower sheet, where $\mathbf{K}_N\approx \kappa_d$, corresponds to the diffusive regime; the narrow transition band between $L^{-}$ and $L^{+}$ is rendered semi-transparently. Panel~(b) shows the corresponding time series $\mathbf{K}_N(t)$ along~$\Gamma$, from which shutdown events appear as brief excursions into the diffusive regime, consistent with crossings of~$L^{-}$. We refer to each such excursion as a \emph{convective shutdown event} and indicate the corresponding minimum of~$\mathbf{K}_N$ by an orange dot in both panels~(a) and~(b). The time series also highlights a clear separation of time scales: $\Gamma$ evolves slowly while $\mathbf{K}_N(t)$ remains on the convective or diffusive sheet, and it undergoes fast transitions when it passes through the narrow intermediate regime between $L^{-}$ and~$L^{+}$.

The periodic solution~$\Gamma$ in Figure~\ref{fig:welander_occ_intro} represents a Welander oscillation in system~\eqref{eq:non_dim_system}. Since $\Gamma$ evolves in a four-dimensional phase space, its projection onto the surface of $\mathbf{K}_N$ may exhibit self-intersections and may cross the loci $L^{-}$ and $L^{+}$ multiple times. We nevertheless identify a convective shutdown event only when $\Gamma$ crosses from the convective sheet, passes through the intermediate regime between $L^{-}$ and $L^{+}$, and enters the diffusive sheet.

\subsection{One-parameter bifurcation diagrams}
\label{section:one_parameter_periodic_solution}
\begin{figure}[t!]
  \centering
  \includegraphics{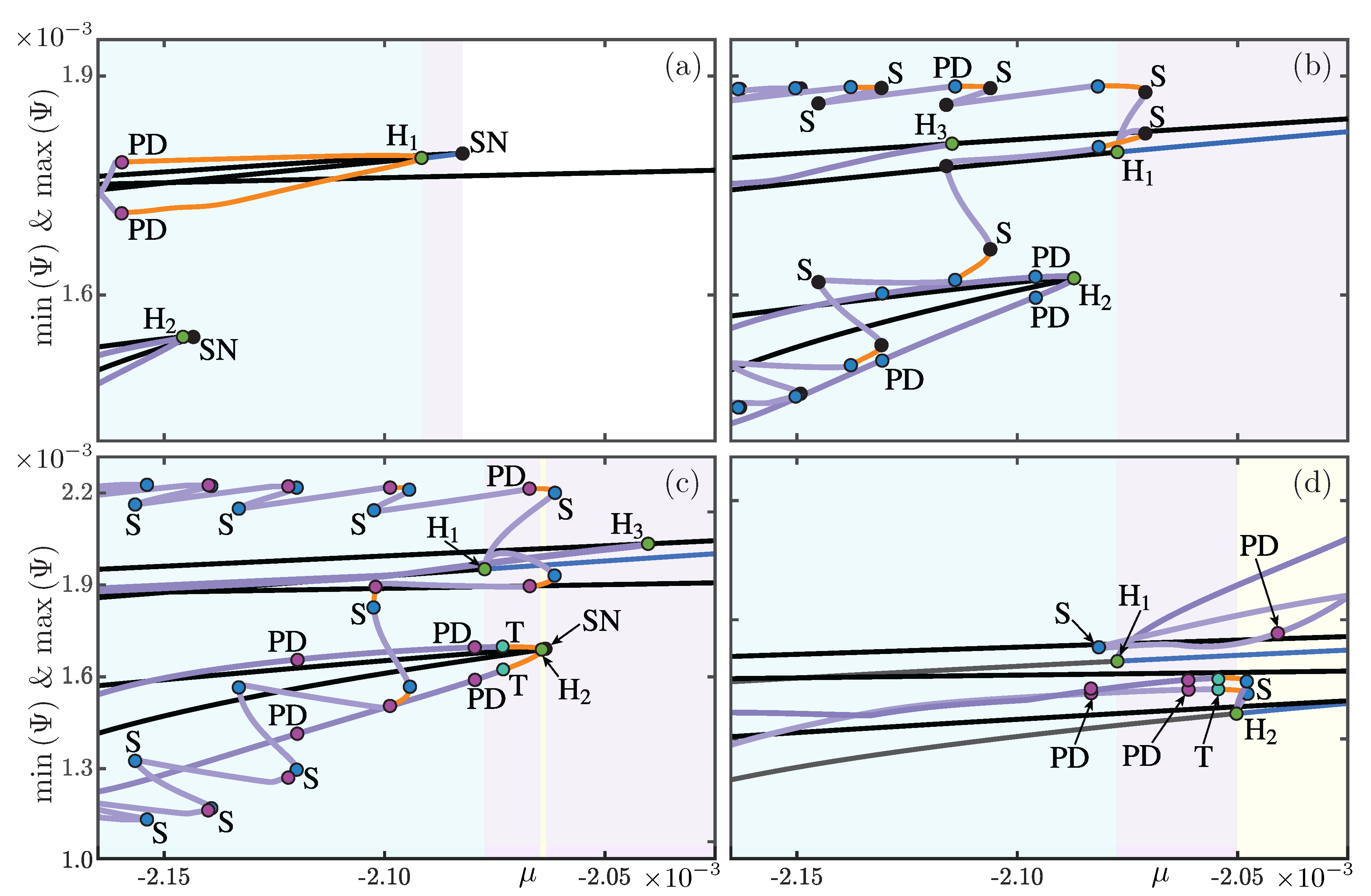}
  \caption{One-parameter bifurcation diagrams of periodic solutions in~$\mu$ for fixed~$\eta$. Stable (orange) and unstable (purple) periodic orbits are shown by plotting extrema of the overturning strength~$\Psi$. All branches bifurcate from Hopf bifurcation points, labeled~$\mathrm{H}_1$, $\mathrm{H}_2$, and $\mathrm{H}_3$ (green dots), and change stability at saddle--node~$\mathrm{S}$ (black dots), period--doubling~$\mathrm{PD}$ (blue dots), and torus~$\mathrm{T}$ (teal dots) bifurcations. Background shading indicates the equilibrium structure: unshaded for a single northward equilibrium, cyan for a single southward equilibrium, purple for coexistence of northward and southward equilibria, and yellow for three equilibria. Panel~(a) is for $\eta=-3.72$, panel~(b) for $\eta=-3.99$, panel~(c) for $\eta=-4.08$, and panel~(d) for $\eta=-4.25$.}
  \label{fig:oneparameter_PO}
\end{figure}

Figure~\ref{fig:oneparameter_PO} shows one-parameter bifurcation diagrams with branches of equilibria and periodic orbits continued in the virtual salinity~$\mu$. Families of periodic solutions are represented by their extrema in overturning strength, $\max(\Psi)$ and $\min(\Psi)$, with stable segments colored in orange and unstable segments in purple. Each family bifurcates from a Hopf bifurcation point~$\mathrm{H}_i$ (green dots, $i=1,2,3$) on an equilibrium branch and subsequently undergoes stability changes at saddle--node bifurcations of periodic orbits~$\mathrm{S}$ (blue dots), period--doubling bifurcations~$\mathrm{PD}$ (purple dots), and torus bifurcations~$\mathrm{T}$ (teal dots). Background shading and equilibrium branch colors are as in Section~\ref{section:onepar_steady_state}, and the $(\mu,\eta)$-range is that of Figure~\ref{fig:EqTwoPar}(c). The southward equilibrium branch from Figure~\ref{fig:EqCurves} lies outside the $\Psi$-range shown and is, therefore, indicated only by cyan shading.

Figure~\ref{fig:oneparameter_PO}(a) for $\eta = -3.72$ shows that, as~$\mu$ decreases, system~\eqref{eq:non_dim_system} transitions from a single northward equilibrium (unshaded), through a bistable regime with a coexisting southward equilibrium (purple shading), to a single southward state (cyan shading). At~$\mathrm{H}_1$, the northward equilibrium loses stability and gives rise to a family of periodic orbits that extends toward lower~$\mu$ and thus appears to lean to the left in the diagram. This family is initially stable but loses stability at a period-doubling bifurcation~$\mathrm{PD}$, beyond which the southward equilibrium is the only attracting state. A second Hopf bifurcation~$\mathrm{H}_2$ generates a family of periodic orbits that remains unstable throughout.

Figure~\ref{fig:oneparameter_PO}(b), for~$\eta = -3.99$, shows bistability between northward and southward equilibria over a range of smaller~$\mu$ (purple shading), while only the southward equilibrium remains stable for more negative~$\mu$ (cyan shading). From~$\mathrm{H}_1$, the northward equilibrium again gives rise to a branch of periodic orbits that extends toward lower values of~$\mu$. This branch is initially unstable, but it alternately gains and loses stability at successive saddle--node~$\mathrm{S}$ and period--doubling~$\mathrm{PD}$ bifurcations, producing a sequence of progressively narrower stable intervals as~$\mu$ decreases. The first of these intervals crosses both the purple and cyan regions, while the subsequent intervals lie entirely within the cyan region. Additional Hopf bifurcation points $\mathrm{H}_2$ and $\mathrm{H}_3$ generate further branches that remain unstable throughout.

In Figure~\ref{fig:oneparameter_PO}(c), for $\eta = -4.08$, a narrow tristable interval (yellow shading) lies between two bistable intervals (purple shading). As in Figure~\ref{fig:oneparameter_PO}(b), the periodic branch that emerges from~$\mathrm{H}_1$ is initially unstable, and its stability alternates at successive saddle--node~$\mathrm{S}$ and period--doubling~$\mathrm{PD}$ bifurcation points. The first stable segment belongs to both purple intervals and the yellow strip in Figure~\ref{fig:oneparameter_PO}(c), while subsequent stable segments lie entirely within the cyan region. The family of periodic solutions emerging from~$\mathrm{H}_2$ is initially stable but quickly loses stability at a torus bifurcation~$\mathrm{T}$. The branch from~$\mathrm{H}_3$ is again unstable throughout.

The final bifurcation diagram in Figure~\ref{fig:oneparameter_PO}(d), for $\eta = -4.25$, shows that as, $\mu$ decreases, the background shading passes successively through yellow, purple, and then cyan intervals. The Hopf bifurcation $\mathrm{H}_1$ remains subcritical; the associated family of periodic orbits is now unstable throughout and extends toward larger values of~$\mu$. In contrast, the family that bifurcates from~$\mathrm{H}_2$ is initially unstable, briefly gains stability at the fold~$\mathrm{S}$ (crossing both the yellow and purple regions), and then loses stability again at the torus bifurcation~$\mathrm{T}$. The Hopf bifurcation~$\mathrm{H}_3$ has moved to smaller $\mu$ and lies outside the plotted range; its unstable branch of periodic orbits enters the shown $\mu$-range of the diagram from the right at $\mu = 0$, turns at the fold~$\mathrm{S}$, and exits again.

Overall, as~$\mu$ decreases, the families of periodic solutions born at the Hopf bifurcation points~$\mathrm{H}_1$, $\mathrm{H}_2$, and~$\mathrm{H}_3$ either remain unstable or admit only a bounded interval of stability. Where stable, these periodic orbits coexist with the southward equilibrium and, for certain~$\eta$, also with the northward equilibrium. Across the cases considered, stronger freshwater forcing (decreasing~$\mu$) weakens the overturning, as is reflected in the plateauing of $\max(\Psi)$ and the pronounced decline of $\min(\Psi)$. We find that only periodic orbits born at~$\mathrm{H}_1$ and~$\mathrm{H}_2$ are stable, whereas the those emerging from~$\mathrm{H}_3$ remain unstable throughout the plotted parameter range.

\subsection{Bifurcation structure of periodic solutions in the $(\mu,\eta)$-plane}
\label{section:bifurcations_single_loop_two_parameter}
\begin{figure}[t!]
  \centering
  \includegraphics{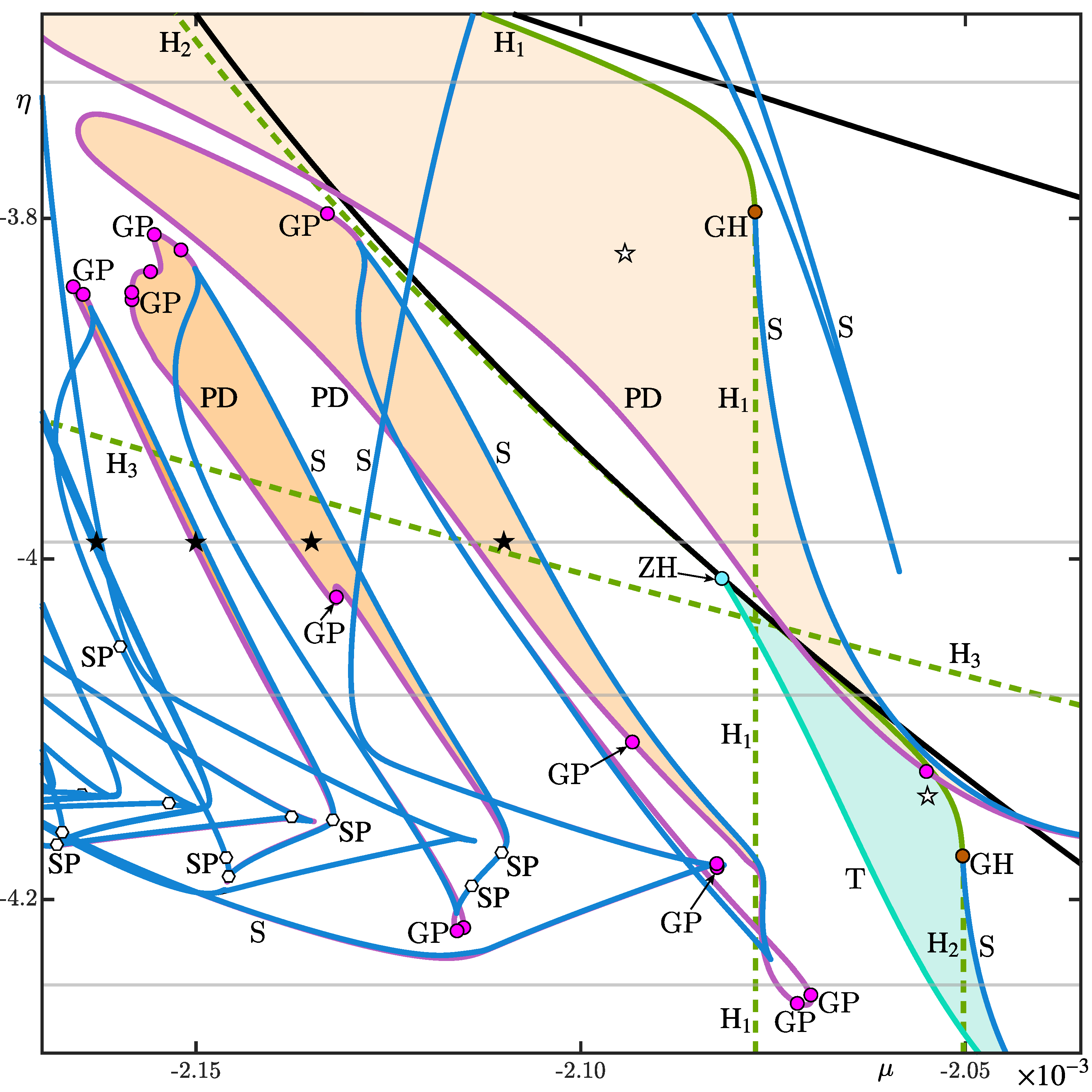}
\caption{\label{fig:two_parameter_po} Two-parameter bifurcation diagram of system~\eqref{eq:non_dim_system} in the $(\mu,\eta)$-plane, over the parameter range in  Figure~\ref{fig:EqTwoPar}(c). Bifurcations of equilibria include Hopf bifurcations~$\mathrm{H_i}$ (green curves; dashed when subcritical and solid when supercritical), and saddle--node bifurcations~$\mathrm{SN}$ (black curves); bifurcations of periodic orbits are saddle--node bifurcations~$\mathrm{S}$ (blue), period-doubling bifurcations~$\mathrm{PD}$ (purple), and torus bifurcations~$\mathrm{T}$ (teal). Codimension-two bifurcations are generalized period-doubling~$\mathrm{GP}$ (magenta dots), saddle-node--period-doubling ~$\mathrm{SP}$ (white marker), generalized Hopf bifurcation~$\mathrm{GH}$ (brown dots), and the zero--Hopf~$\mathrm{ZH}$ (cyan dot). Orange--yellow shading indicates oscillations with predominantly active North Atlantic convection, and teal shading oscillations with tropical convection. Grey lines mark the $\eta$-values used in Figure~\ref{fig:oneparameter_PO}; white stars indicate the parameter values of the representative solutions in Figure~\ref{fig:rep_one_orbits_start}, and black stars those in Figure~\ref{fig:seq_periodic_solutions}.}
\end{figure}

Continuation of the bifurcation points~$\mathrm{H}_i$ and~$\mathrm{SN}$ bifurcations of equilibria, as well as the saddle--node~$\mathrm{S}$, period-doubling~$\mathrm{PD}$, and torus~$\mathrm{T}$ bifurcations of periodic orbits results in the bifurcation diagram shown in Figure~\ref{fig:two_parameter_po}. Here, the color scheme of the curves in the $(\mu,\eta)$-plane is as in Figures~\ref{fig:EqCurves} and~\ref{fig:EqTwoPar}, with the Hopf curves~$\mathrm{H}_i$ drawn dashed where they are subcritical and solid where they are supercritical; we do not distinguish the criticality of the period-doubling curves~$\mathrm{PD}$.

Certain points of codimension-two bifurcation points organize the structure of bifurcation diagram. At the generalized Hopf bifurcation points~$\mathrm{GH}$ on the Hopf curves~$\mathrm{H}_1$ and~$\mathrm{H}_2$, the Hopf bifurcation changes criticality and a curve~$\mathrm{S}$ of saddle--node bifurcations of periodic orbits emerges. We also find additional curves~$\mathrm{S}$ that are not connected to a point $\mathrm{GH}$ in the parameter range of Figure~\ref{fig:two_parameter_po}. Generalized period-doubling points~$\mathrm{GP}$ on period-doubling curves~$\mathrm{PD}$ can likewise be the origin of additional saddle--node bifurcation curves of the period-doubled (two-loop) periodic orbits~\cite{kuznetsov1998elements}; these curves are not shown to keep Figure~\ref{fig:two_parameter_po} simple. A tangency between the curves~$\mathrm{S}$ and~$\mathrm{PD}$ takes place at the codimension-two points~$\mathrm{SP}$ of simultaneous saddle-node and period doubling bifurcation, which occur predominantly for smaller~$\eta$. Finally, a zero--Hopf (or saddle-node Hopf) bifurcation point~$\mathrm{ZH}$ occurs at a tangential intersection of the Hopf curve~$\mathrm{H}_2$ with a saddle--node curve~$\mathrm{SN}$. It is of type~$\rm{IV}$ in the notation of~\cite{kuznetsov1998elements}, changes the criticality of~$\mathrm{H}_2$, and gives rise to the (supercritical) torus bifurcation curve~$\mathrm{T}$.

The bifurcation curves delineate regions of the $(\mu,\eta)$-plane in which system~\eqref{eq:non_dim_system} supports stable oscillations. Specifically, the torus bifurcation curve~$\mathrm{T}$, together with a saddle--node curve~$\mathrm{S}$ and the supercritical segment of~$\mathrm{H}_2$, bounds a region (shaded teal) that admits a stable periodic orbit. There are also five additional regions (shaded yellow) of stable periodic solutions, one of which is very narrow; they are each bounded by segments of a saddle--node curve~$\mathrm{S}$ and a period-doubling curve~$\mathrm{PD}$. These yellow regions become progressively smaller and accumulate as $\mu$ decreases, a structure organized by homoclinic tangencies associated with a codimension-two Shilnikov--Hopf point~\cite{kuznetsov1998elements} at $(\mu,\eta) \approx (-2.216\times10^{-3}, -3.8552)$, which lies just outside the parameter range of Figure~\ref{fig:two_parameter_po}. The accumulation of period-doubling bifurcations near this point suggests a route to chaotic dynamics, which we investigate in Section~\ref{section:chaotic_pulse_welander}.

\begin{figure}[t!]
  \centering
  \includegraphics{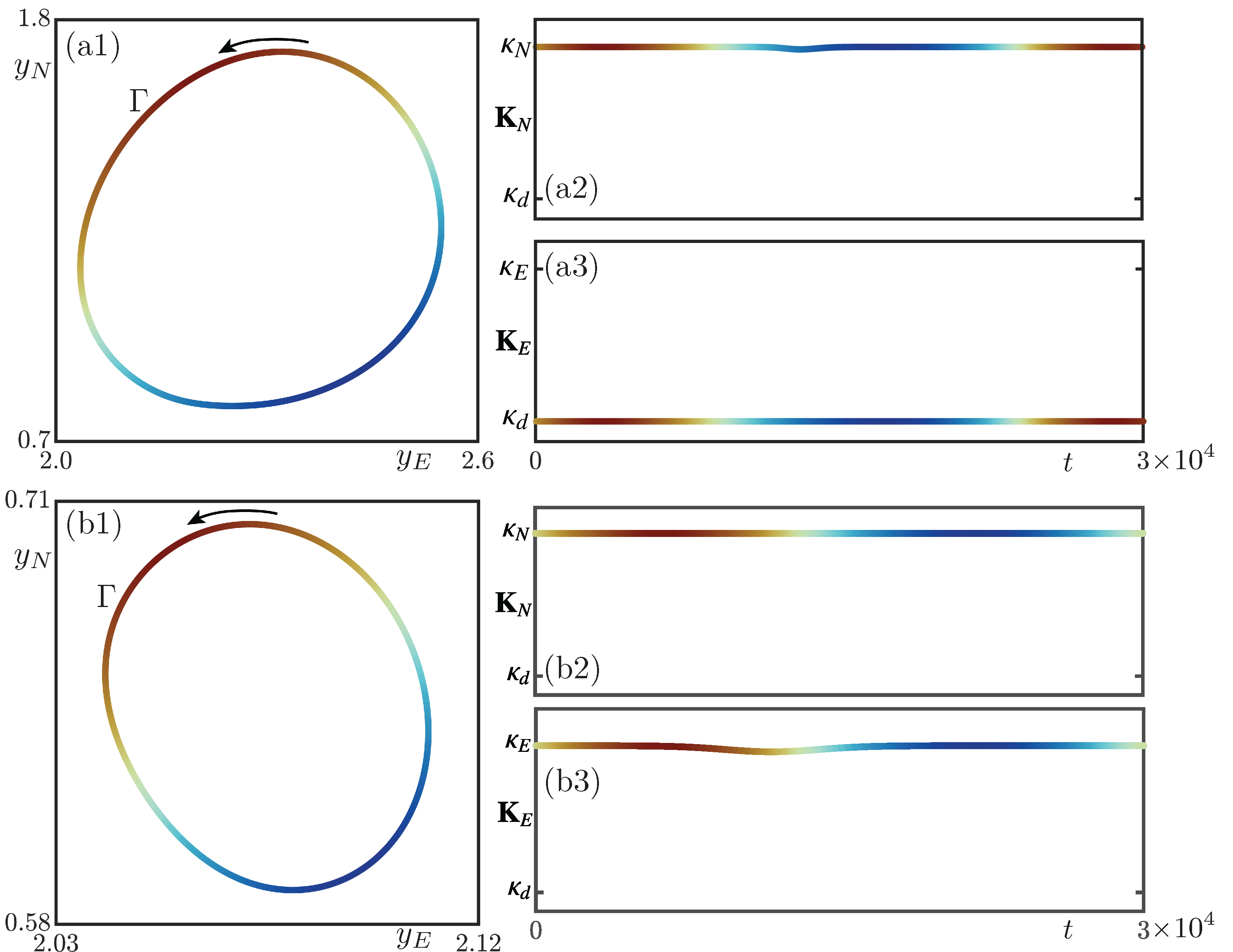}
  \caption{\label{fig:rep_one_orbits_start}
  Two types of non-Welander oscillations, shown over one period and colored by the overturning strength~$\Psi$ from weak (blue) to strong (red). Panel~(a1), for $(\mu,\eta)=(-2.0978\times10^{-3},-3.82)$, shows the projection of~$\Gamma$ onto the $(y_N,y_E)$-plane; panel~(a2) shows the time series of~$\mathbf{K}_N(t)$; and panel~(a3) shows the time series of~$\mathbf{K}_E(t)$. Panels~(b1)--(b3) are shown in the same manner for $(\mu,\eta)=(-2.0978\times10^{-3},-3.82)$.}
  \end{figure}

\subsection{Representative stable oscillations}
\label{section:representations_welander_orbits}
We now illustrate representative periodic solutions of system~\eqref{eq:non_dim_system} within the regions identified in Figure~\ref{fig:two_parameter_po}, at the parameter values indicated by white and black star markers. Each periodic solution~$\Gamma$ is shown over one period by its projection onto the $(y_N,y_E)$-plane and by the time series of~$\mathbf{K}_N(t)$ and either~$\mathbf{K}_E(t)$ or~$\Psi(t)$. All representations are coloured by the overturning strength~$|\Psi|$, from weak (dark blue) to strong (dark red), and convective shutdown events are marked by orange dots. We distinguish two types of periodic solutions: oscillations without convective shutdown events, and Welander oscillations with convective shutdown events. 
\begin{figure}[t!] 
  \centering
    \includegraphics{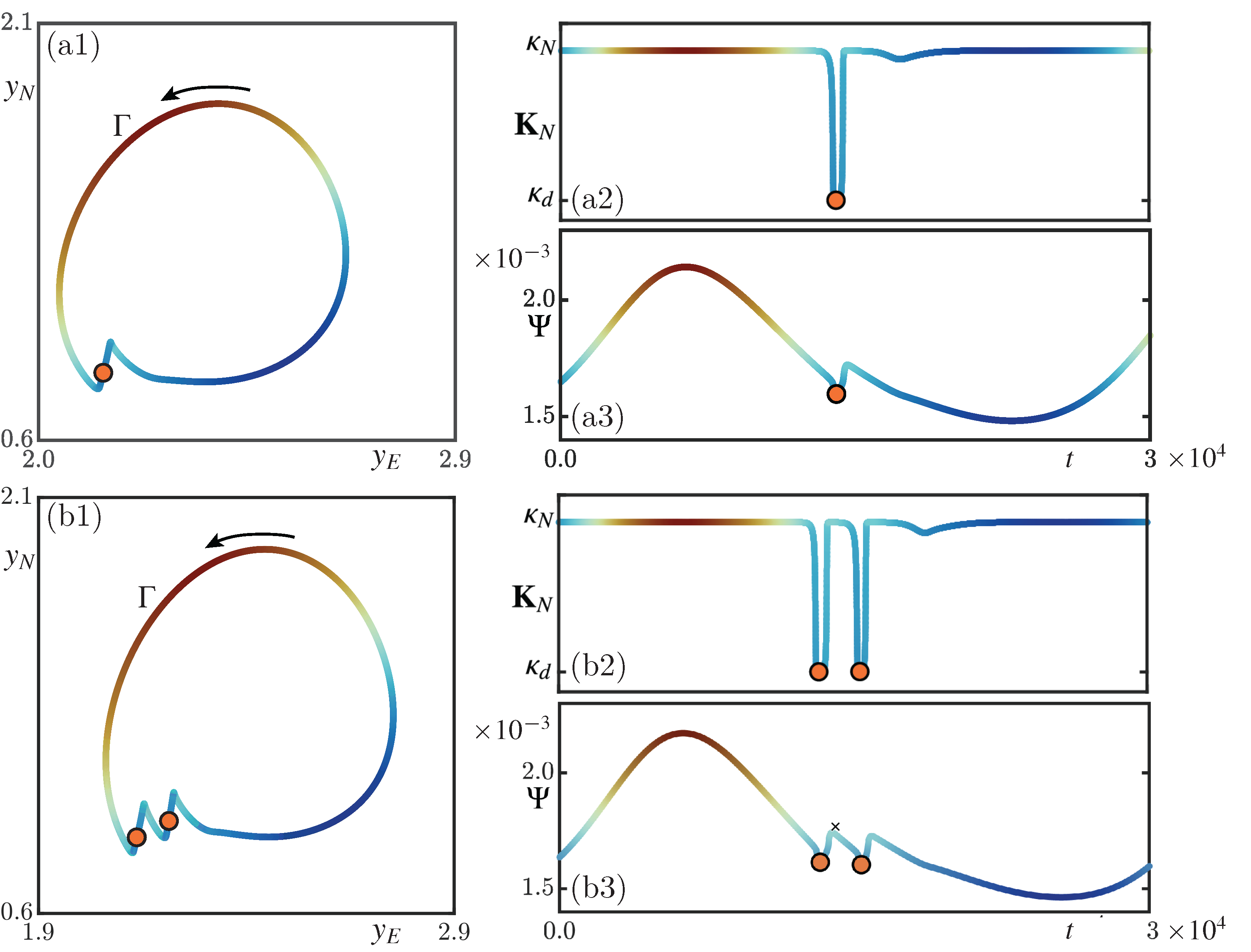}
  \caption{\label{fig:seq_periodic_solutions}
Four Welander oscillations over one period, represented as in Figure~\ref{fig:rep_one_orbits_start} but with the time series of~$\Psi(t)$ replacing that of~$\mathbf{K}_E(t)$; convective shutdown events are indicated by orange dots. Panels~(a1)--(a3) show a single shutdown event for $(\mu,\eta)=(-2.1086\times10^{-3},-3.99)$. Panels~(b1)--(b3) show two shutdown events for $(\mu,\eta)=(-2.135\times10^{-3},-3.99)$. Panels~(c1)--(c3) show three shutdown events for $(\mu,\eta)=(-2.15\times10^{-3},-3.99)$. Panels~(d1)--(d3) show four shutdown events for $(\mu,\eta)=(-2.163\times10^{-3},-3.99)$.
}
\end{figure}


\begin{figure}[t!]\ContinuedFloat
  \centering
    \includegraphics{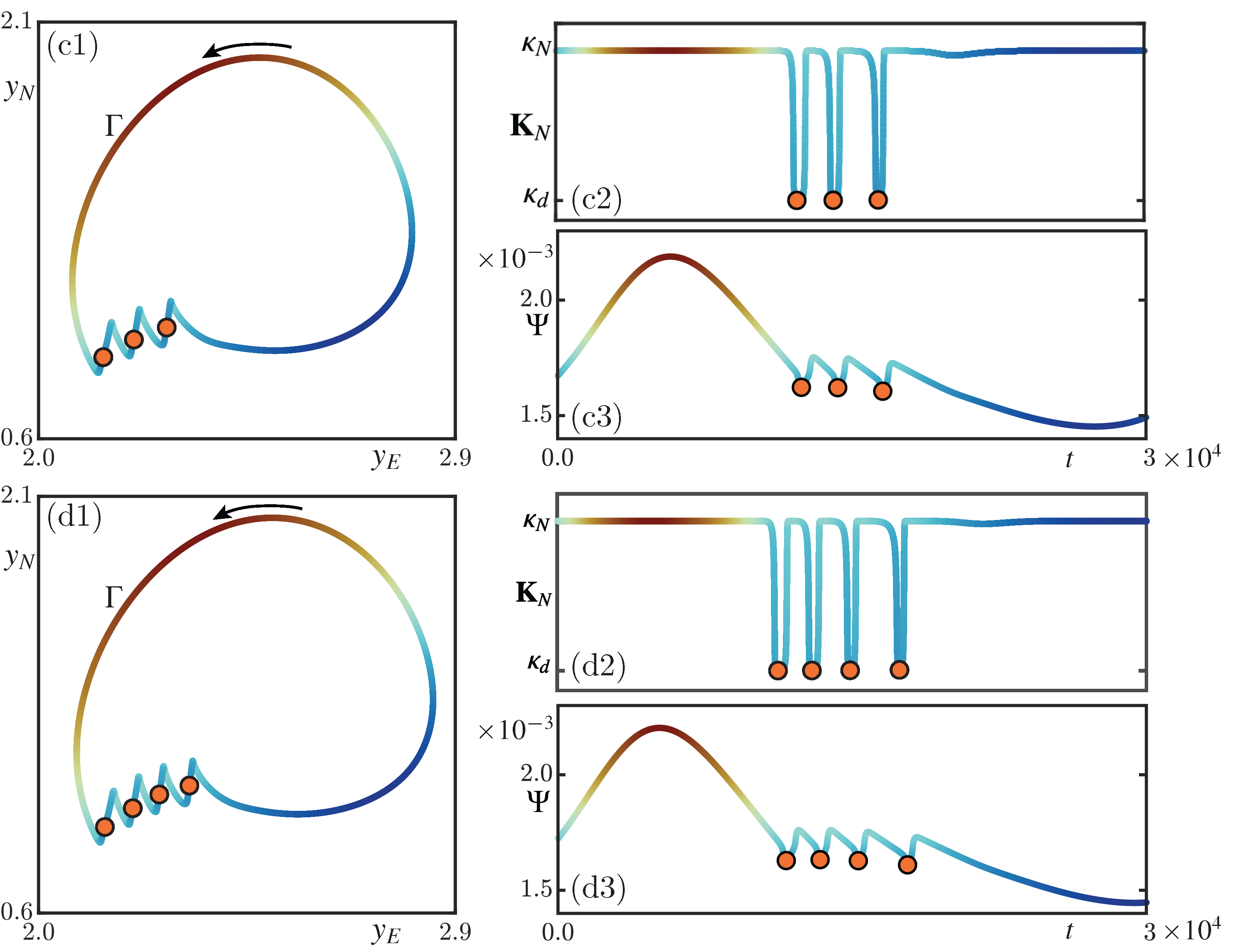}
    \caption{\label{fig:seq_periodic_solutions_cont}
    continued. }
\end{figure}

Figure~\ref{fig:rep_one_orbits_start} shows representative periodic solutions at the white stars in the teal region and in the rightmost yellow region. Both solutions have northward advection~($\Psi > 0$), but neither exhibits convective shutdown events. Panels~(a1) and~(b1) of Figure~\ref{fig:rep_one_orbits_start} reveal that these two solutions are indistinguishable in the~$(y_N,y_E)$-plane, yet their physical properties differ. In the rightmost yellow region, panels~(a2)--(a3) show that North Atlantic convection remains active throughout the cycle, with~$\mathbf{K}_N(t)$ near the convective level~$\kappa_N$, while tropical mixing stays at the diffusive level~$\kappa_d$. Along~$\Gamma$ in the teal region, on the other hand, panels~(b2)--(b3) show that both the North Atlantic and tropical convection remain active. In neither case does $\mathbf{K}_N(t)$ or~$\mathbf{K}_E(t)$ switch between the convective and diffusive levels and, hence, neither periodic solution constitutes a Welander oscillation.

In the additional yellow regions in Figure~\ref{fig:two_parameter_po}, however, one finds four distinct Welander oscillations; these are shown in Figure~\ref{fig:seq_periodic_solutions} at the parameter values indicated by the black stars in Figure~\ref{fig:two_parameter_po}. The simplest case, shown in Figure~\ref{fig:seq_periodic_solutions}(a1)--(a3) for $(\mu,\eta) = (-2.1086\times 10^{-3}, -3.99)$, is the oscillation already displayed in Figure~\ref{fig:welander_occ_intro}; however it is now shown in terms of $\Psi$ and coloured by the overturning strength~$|\Psi|$. Its single shutdown event appears as a brief zig-zag in the $(y_N,y_E)$-projection at relatively small values of $y_N$ and $y_E$, which corresponds to a rapid drop of $\mathbf{K}_N(t)$ to the diffusive level~$\kappa_d$ followed by a fast recovery toward~$\kappa_N$. In~$\Psi(t)$ the convective shutdown event manifests itself as a short-lived dip that occupies only a small fraction of the period; this dip has features of slow--fast dynamics reminiscent of those found in the Welander model~\cite{bailie2025detailed}.

The remaining three oscillations involve progressively more shutdown events per period. Panels~(b1)--(b3) of Figure~\ref{fig:seq_periodic_solutions} show the case with two shutdown events: here, $\Gamma$ exhibits two zig-zags in panel~(b1), with corresponding short dips in~$\mathbf{K}_N(t)$ in panel~(b2) and a pair of minima in~$\Psi(t)$ in panel~(b3). Panels~(c1)--(c3) and~(d1)--(d3) show the analogous behaviour for three and four shutdown events, respectively. Across the $(\mu,\eta)$-plane in Figure~\ref{fig:two_parameter_po}, as the virtual salinity $\mu$ becomes more negative (an increase in freshwater influx), each successive yellow region corresponds to one additional shutdown event per period. Individual shutdown--recovery episodes remain brief but, as their number increases, a progressively larger fraction of the cycle is spent at intermediate overturning strengths: a single episode occupies about $5$\% of the period in Figure~\ref{fig:seq_periodic_solutions}(a3), whereas four shutdown--recovery episodes occupy about $25$\% in panel~(d3).

\begin{figure}[t!]
  \centering
  \includegraphics{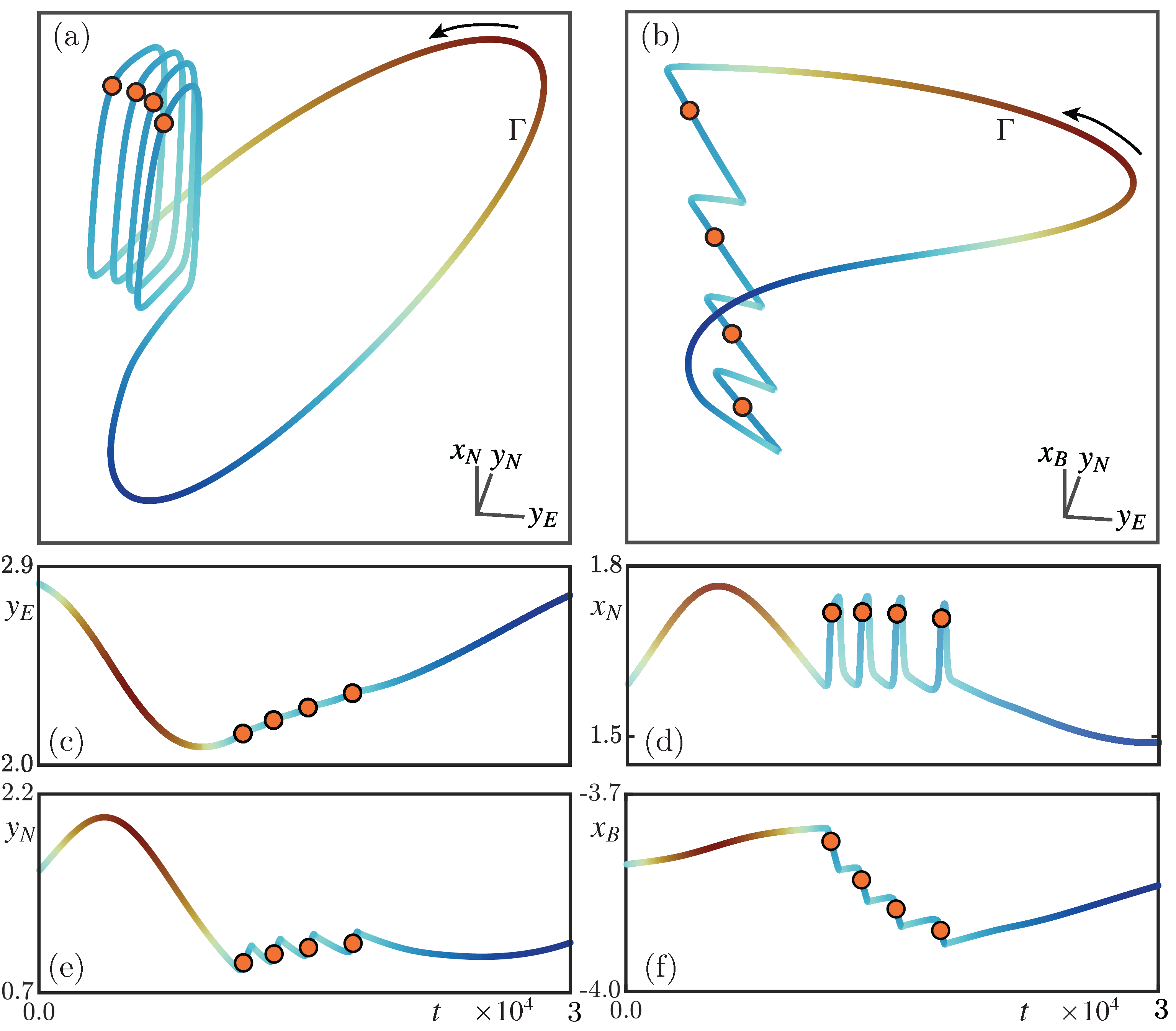}
\caption{\label{fig:3D_po_representations} Welander oscillation~$\Gamma$ for $(\mu,\eta)=(-2.163\times10^{-3},-3.99)$ from Figure~\ref{fig:seq_periodic_solutions}. Panel~(a) shows the projection onto $(y_N,y_E,x_N)$-space and panel~(b) onto $(y_N,y_E,x_B)$-space; panels~(c)--(f) are time series of the state variables $(y_E,x_N,y_N,x_B)$.}
\end{figure}

From a physical perspective, each Welander oscillation in system~\eqref{eq:non_dim_system} consists of two distinct phases: a slow preconditioning phase and a fast shutdown--recovery phase. This is illustrated in Figure~\ref{fig:3D_po_representations} for the Welander oscillation with four convective-shutdown events from Figure~\ref{fig:seq_periodic_solutions}(d1)--(d3). Projections onto the $(y_N,y_E,x_N)$- and $(y_N,y_E,x_B)$-spaces in in Figure~\ref{fig:3D_po_representations}, together with time series of each variable, reveal the physical mechanisms behind these two phases. During the slow preconditioning phase, starting at peak overturning (as indicated by dark red), the subpolar salinity~$y_N$, the tropical salinity~$y_E$, and the North Atlantic temperature~$x_N$ all decrease, while the deep-water temperature~$x_B$ increases very gradually. This combination narrows the density contrast between the surface and deep-water boxes until it crosses the convective-adjustment threshold~$\eta$, at which point the fast shutdown--recovery phase is initiated. During this phase, $x_N$ rises sharply while $x_B$ drops; this is visible in the $(y_N,y_E,x_N)$-projection in panel~(a) as a loop. Subsequently, $x_N$ decreases quickly, but $x_B$ evolves on a much longer timescale and recovers only partially, producing the step-like descent in the $(y_N,y_E,x_B)$-projection in panel~(b). The convective shutdown event itself occurs just before the maximum of~$x_N$ and halfway down the sharp drop of~$x_B$. Notably, each incomplete recovery of~$x_B$ leaves the density contrast close to the threshold~$\eta$, so that the next shutdown event is triggered more readily; this causes multiple shutdown events to cluster along~$\Gamma$. The sequence of shutdown events terminates only once~$y_N$ and~$y_E$, shown in panels~(c) and~(e), have evolved sufficiently to restore the density contrast above~$\eta$, and the slow preconditioning phase resumes.

The deep temperature~$x_B$ thus acts as a slowly evolving memory variable in system~\eqref{eq:non_dim_system}: it couples successive shutdown events through its incomplete recovery, as in Rahmstorf-type box models~\cite{rahmstorf1996freshwater}. The resulting pattern of slow preconditioning followed by rapid collapse of convection is reminiscent of view the so-called Dansgaard--Oeschger variability~\cite{rahmstorf2002ocean}, albeit in a highly idealized setting that omits sea-ice feedbacks, long-range teleconnections, and interhemispheric coupling~\cite{dokken2013dansgaard}.

\begin{figure}[ht!]
  \centering
  \includegraphics[scale=1.0]{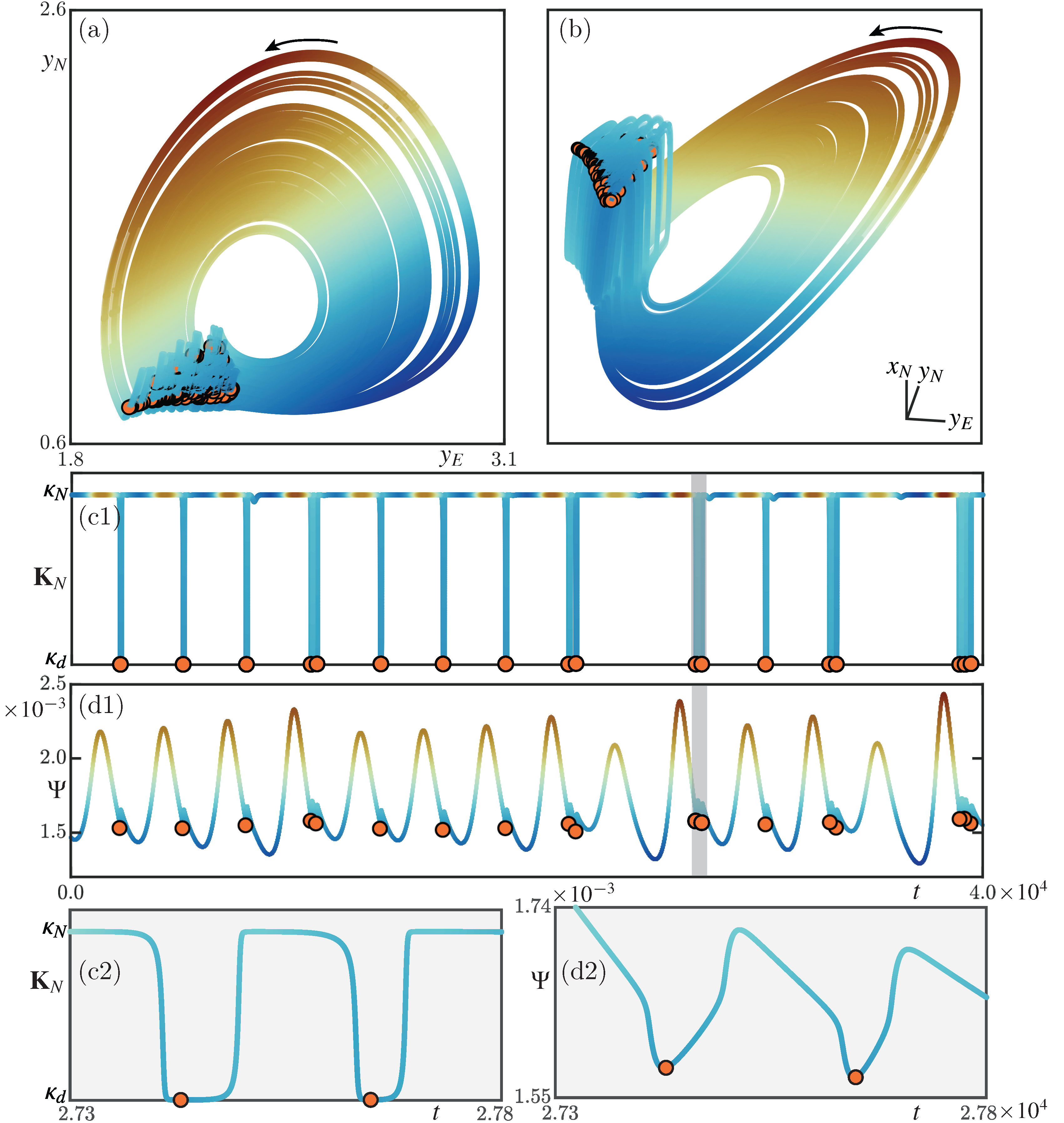}
\caption{\label{fig:rep_chaotic_left} Chaotic Welander oscillation for $(\mu,\eta) = (-2.128\times10^{-3}, -3.99)$, with convective shutdown events marked by orange dots. Panel~(a) shows the projection onto the $(y_N,y_E)$-plane, colored by the overturning strength~$\Psi$; panel~(b) onto $(y_N,y_E,x_N)$-space; panels~(c1) and~(d1) are time series of~$\mathbf{K}_N(t)$ and~$\Psi(t)$; and panel~(c2) and~(c2) are enlargements for $t\in[0,4]\times10^{4}$. The grey band, in panels~(c1) and~(d1), indicates the time interval magnified in panels~(c2) and~(d2), which show a single group of two consecutive shutdown events.}
\end{figure}

\section{Chaotic Welander oscillations}
\label{section:chaotic_pulse_welander}
The convective shutdown events considered thus far recur at regular intervals along stable periodic orbits of system~\eqref{eq:non_dim_system}. As noted in Section~\ref{section:bifurcations_single_loop_two_parameter}, the accumulation of period-doubling bifurcations~$\mathrm{PD}$ provides a mechanism for chaotic dynamics. We now show that shutdown events can occur irregularly; we refer to this type of dynamics as chaotic Welander oscillations, which are associated with chaotic attractors of system~\eqref{eq:non_dim_system}.

Figure~\ref{fig:rep_chaotic_left} shows an example of a trajectory on the chaotic attractor that exists for $(\mu,\eta) = (-2.128\times10^{-3}, -3.99)$. Panel~(a) displays its projection onto the~$(y_N,y_E)$-plane, panel~(b) shows the trajectory in~$(y_N,y_E,x_N)$-space, and panels~(c1) and~(d1) show the corresponding time series of~$\mathbf{K}_N(t)$ and~$\Psi(t)$ over~$t\in[0,4 \times10^{4}]$; panels~(c2) and~(d2) are enlargements around a single cluster of shutdown events. The trajectory in Figure~\ref{fig:rep_chaotic_left}(a) and~(b) spirals anti-clockwise, and appears to accumulate on a thin sheet, from which it is intermittently ejected during convective shutdown events (orange dots) before returning to the sheet. As in the periodic Welander oscillation shown in Figure~\ref{fig:3D_po_representations}, shutdown events cluster at smaller values of~$y_N$ and~$y_E$; in contrast to the periodic case, however, infinitely many such events occur along the chaotic trajectory. The time series of~$\mathbf{K}_N$ in Figure~\ref{fig:rep_chaotic_left}(c1) forms a `barcode' of irregular shutdown events, where each dip to the diffusive level~$\kappa_d$ coincides with a local minimum of~$\Psi$ in panel~(d1). Panels~(c2) and~(d2) enlarge a cluster with two shutdown events, each of which exhibits the same sawtooth structure as for the periodic case. The picture that emerges is that larger local maxima of~$\Psi$ precede larger clusters of shutdown events, whereas smaller local maxima are followed not by shutdown events but by a deeper overall dip in~$\Psi$. Overall, the chaotic attractor switches unpredictably between clusters containing one, two, or three closely spaced events.

\begin{figure}[t!]
  \centering
  \includegraphics{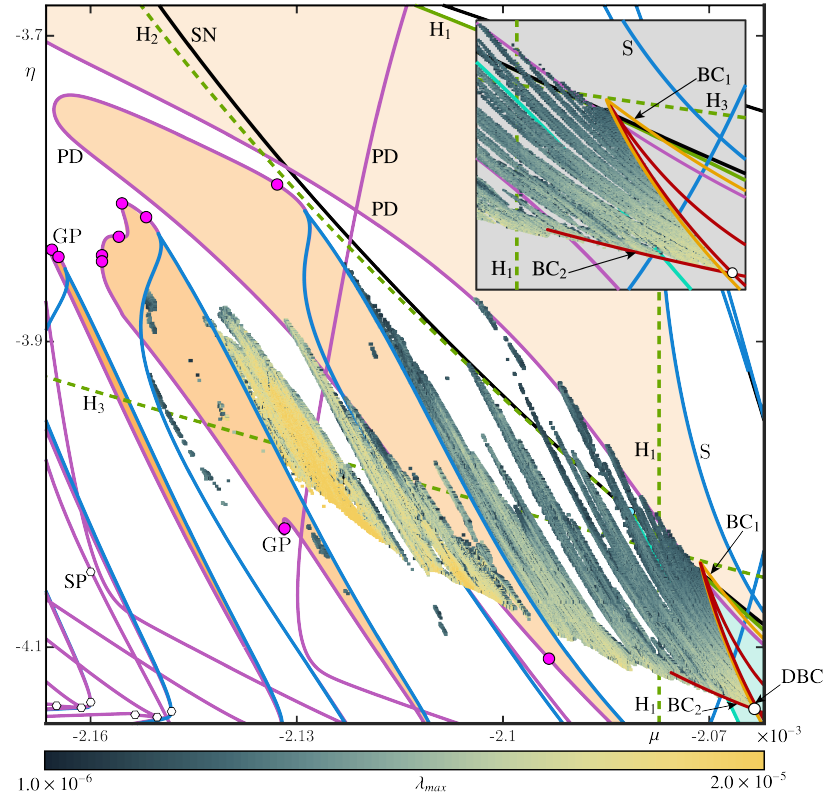}
\caption{\label{fig:two_parameter_chaos} Two-parameter bifurcation diagram in the $(\mu,\eta)$-plane from Figure~\ref{fig:two_parameter_po} over a slightly enlarged parameter range, with superimposed values of the maximal Lyapunov exponent~$\lambda_{\max}$ as given by the color bar. Also shown are curves~$\mathrm{BC}_1$ (yellow) and~$\mathrm{BC}_2$ (red), which meet at the point~$\mathrm{DBC}$ (white dot), as illustrated by the inset. }
\end{figure}

Having identified chaotic Welander oscillations at a single parameter point, we now ask where they exist in the~$(\mu,\eta)$-plane. To this end, we compute the largest Lyapunov exponent~$\lambda_{\max}$ over an~$800\times800$ grid with a standard variational algorithm~\cite{benettin1976kolmogorov} as implemented in \texttt{ChaosTools.jl}~\cite{Julia_2017,DatserisParlitz2022}. Figure~\ref{fig:two_parameter_chaos} shows the resulting map, where shading indicates parameter values with~$\lambda_{\max}$ above~$10^{-6}$; values below this threshold are treated as zero to account for numerical noise. The set of detected chaotic dynamics is quite large, and the chaotic attractor shown in Figure~\ref{fig:rep_chaotic_left} is indicated by the black star in Figure~\ref{fig:two_parameter_chaos}.

The parameter set with chaotic dynamics is interspersed with periodic windows in which trajectories converge to stable periodic orbits of higher period. These windows are bounded by saddle--node and period-doubling bifurcation curves (not shown) of the corresponding periodic orbits~\cite{kuznetsov1998elements}; here, they are identified as ranges of non-positive values of~$\lambda_{\max}$. The set of positive $\lambda_{\max}$ overlaps parts of the second and third orange--yellow regions, so that periodic and chaotic Welander oscillations coexist; it also extends into the teal region, where non-Welander oscillations are found.

For parameter values in the lower-right corner of Figure~\ref{fig:two_parameter_chaos}, we find that the basin boundary of the chaotic attractor is formed by the stable manifold of a saddle periodic orbit. The attractor is destroyed in a boundary crisis when it collides with this basin boundary, which occurs here at a homoclinic tangency between the stable and unstable manifolds of the periodic orbit~\cite{grebogi1983crises,grebogi1987chaos}. We compute and continue two such tangencies with an implementation of Lin's method~\cite{krauskopf2008lin}, which results in the curves~$\mathrm{BC}_1$ and~$\mathrm{BC}_2$. The set of positive~$\lambda_{\max}$ extends from the codimension-two double-boundary-crisis point~$\mathrm{DBC}$, where~$\mathrm{BC}_1$ and~$\mathrm{BC}_2$ meet; see the inset of Figure~\ref{fig:two_parameter_chaos}. The two curves bound only part of the parameter set with chaotic dynamics. Its remaining lower edge appears to be organized by additional curves, which suggests the presence of further boundary crises. A detailed analysis of these bifurcations will require the identification of additional saddle periodic orbits and the continuation of connecting orbits between them, which is beyond the scope of the present paper.

\section{Conclusion}
\label{section:conclusion}
We introduced a conceptual AMOC model that couples Stommel-type density-driven overturning with Welander-style convective adjustment. The model~\eqref{eq:non_dim_system} takes the form of a four-dimensional ordinary differential equation and thus provides a setting in which steady-state overturning and Welander oscillations can be analyzed together. Specifically, we performed a comprehensive one- and two-parameter bifurcation analysis in the virtual salinity flux~$\mu$ and the density threshold~$\eta$. This revealed how different stable solutions, and the bifurcations that separate them, are organized in the~$(\mu,\eta)$-plane. In this way, the model clarifies how changes in freshwater forcing can weaken overturning, induce repeated convective shutdown, and ultimately generate irregular AMOC variability.

As~$\mu$ decreases, system~\eqref{eq:non_dim_system} passes from a regime with a single stable equilibrium and strong northward overturning through several intervals of multi-stability with coexisting equilibria, after which only a single southward-overturning equilibrium remains stable for sufficiently negative~$\mu$. The classical Welander mechanism examined in~\cite{bailie2025detailed} carries over to self-sustained oscillations that exhibit a clear separation of timescales: the overturning strength~$\Psi$ evolves on a slow millennial timescale, while convective shutdown events occur on decadal-to-centennial timescales. A feature of both periodic and chaotic Welander oscillations is the incomplete recovery of the deep temperature after each shutdown, which leaves the density contrast close to the threshold~$\eta$ and thereby promotes further shutdown events. In the periodic case, the number of shutdown events per cycle increases as the freshwater forcing strengthens (lower virtual salinity flux~$\mu$). In the parameter range considered here, periodic orbits may exhibit up to four such events, so that approximately~$25\%$ of the cycle is spent in shutdown--recovery episodes at intermediate overturning strengths. The model also supports chaotic Welander oscillations with irregularly timed shutdown events, which may, for some values of~$\mu$ and~$\eta$, coexist with a periodic Welander oscillation.

The periodic Welander oscillations are organized by bounded regions of stable periodic orbits in the $(\mu,\eta)$-plane, each delimited by saddle--node and period-doubling bifurcation curves that accumulate for lower values of~$\mu$. Within the parameter range considered, chaotic Welander oscillations occupy a substantial set in the $(\mu,\eta)$-plane that is bounded, in part, by curves of boundary crises bifurcations. Overall, our results show that a conceptual climate model such as system~\eqref{eq:non_dim_system} can support both regular and irregular convective shutdown dynamics within a purely deterministic setting.

The non-dimensional parameter region explored here corresponds to a physically comparatively large range (of~$(F_N,\Delta\rho)\in[-2.77,0]~\mathrm{Sv}\times[-4.55,0]~\mathrm{kg\,m^{-3}}$;  see Section~\ref{section:model_reduction}), as a result of the basin-scale box volumes and the choice of~$\varepsilon$. In fact, the precise locations of bifurcation curves in the~$(\mu,\eta)$-plane are sensitive to the smoothing parameter~$\varepsilon$ and, hence, our analysis identifies the qualitative bifurcation structure of system~\eqref{eq:non_dim_system}, rather than physically definitive threshold values; it should be interpreted in the context of the chosen regularization and basin volumes.

The analysis presented here can be extended in several natural directions. One is to subdivide the subpolar North Atlantic into distinct ocean basins, such as the Irminger and Labrador seas, each with its own convective-adjustment process. This would make it possible to examine the on--off competition between different convection sites~\cite{neff2023bifurcation} in the presence of advective feedbacks; such a model extension may give rise to additional oscillatory modes involving synchronous and asynchronous shutdown dynamics~\cite{kuznetsov1998elements}. A second direction is to include both hemispheres by adding a South Atlantic box, thereby enabling studies of interhemispheric density gradients and their influence on overturning strength~\cite{castellana2019transition}. Finally, it would be of interest to incorporate stochastic atmospheric forcing to investigate noise-induced transitions between the coexisting periodic and chaotic shutdown regimes identified in system~\eqref{eq:non_dim_system}.

\section*{Acknowledgments}
The authors thank Andrew Keane for helpful discussions during the development of the model. This research was supported in part by the Royal Society Te~Ap\={a}rangi Marsden Fund under grant No.~19-UOA-223.

\bibliographystyle{plain}
\bibliography{BSK_references}
\end{document}